\def\bu{\bullet}
\def\marker{\>\hbox{${\vcenter{\vbox{
					\hrule height 0.4pt\hbox{\vrule width 0.4pt height 6pt
						\kern6pt\vrule width 0.4pt}\hrule height 0.4pt}}}$}\>}

\documentclass[12pt]{article}
\usepackage{indentfirst}
\usepackage{epsfig}
\usepackage{fullpage}
\usepackage{amsmath,amsfonts,amssymb,amsthm,graphics,amscd}
\usepackage{pgf,tikz}
\usepackage{graphicx}
\usepackage{subfigure}
\usepackage{amsmath,amsfonts,amssymb,amsthm,graphics,amscd,pst-plot,pstricks}
\usepackage{graphicx}
\usepackage{float}

\newtheorem {Theorem}  {Theorem}
\newtheorem {Lemma}{Lemma}[section]

\theoremstyle{definition}
\newtheorem{Definition}{Definition}

\newtheorem{Conjecture}{Conjecture}

\usepackage[labelfont=bf]{caption}
\newcommand{\D}{\Delta}

\newcommand{\phiv}{\varphi}

\newcommand{\CC}{\mathcal{C}}

\usepackage[
pdfauthor={},
pdftitle={An  improved bound on the Overfull Conjecture for graphs with small minimum degree and large maximum degree},
pdfstartview=XYZ,
bookmarks=true,
colorlinks=true,
linkcolor=blue,
urlcolor=blue,
citecolor=blue,
bookmarks=false,
linktocpage=true,
hyperindex=true
]{hyperref}

\usepackage[mathlines]{lineno}
\begin{document}
%\pagewiselinenumbers
	
    \begin{center}	
		{\Large \bf An improved bound on the Overfull Conjecture for graphs with small minimum degree and large maximum degree}
		\vspace{6mm}
		
		Xuli Qi$^*$,\,\,Chunhui Ge,\,\,Yanrui Feng
        
		\vspace{3mm}
		\baselineskip=0.25in
		{\it Department of Mathematics, Hebei University of Science and Technology,}  \\
		{\it Shijiazhuang 050018, P.R. China}\\
		\vspace{3mm}
	\end{center}	
    \footnotetext{Corresponding author. E-mail address: xuliqi@hebust.edu.cn}
	\date{}
	
	\begin{abstract}
		Let $G$ be a simple graph with order $n$, maximum degree $\D(G)$, minimum degree $\delta(G)$, and chromatic index $\chi'(G)$. A graph $G$ is called {\em $\D$-critical} if $\chi'(G)=\D(G)+1$ and $\chi'(H) < \chi'(G)$ for every proper subgraph $H$ of $G$, and $G$ is overfull if $\left|E(G)\right|>\Delta(G)\lfloor \frac{n}{2}\rfloor$. In 1986, Chetwynd and Hilton proposed the Overfull Conjecture, which is equivalent to the following $\D$-critical form: Every $\D$-critical graph $G$ with $\D(G) > \frac{n}{3}$ is overfull. The Overfull Conjecture has many implications. For example, it implies a polynomial-time algorithm for determining the chromatic index of graphs $G$ with $\D(G) > \frac{n}{3}$, and implies several longstanding conjectures in the area of graph edge coloring.
		Recently, Cao, Chen, Jing, and Shan (SIAM J. Discrete Math. 2022) verified the Overfull Conjecture for $\D$-critical graphs $G$ when $\D(G)-\frac{7\delta(G)}{4}\ge \frac{3n-17}{4}$. In this paper, we establish the required bound  $\D(G)-\frac{5\delta(G)}{3}\ge\frac{2n-7}{3}$.
		\vspace{2mm}	
	\end{abstract}

	{\tiny {\normalsize {\Large }}}	\section{Introduction}
	We consider simple and connected graphs in this paper. %Some definitions and results on simple graphs are mentioned below.
	We generally follow the book \cite{SSTF2012} of Stiebitz et al. for notation and terminology.
	Let $G=(V(G),E(G))$ be a simple graph, where $V(G)$ is the vertex set and $E(G)$ is the edge set of $G$.
	Denote by $\Delta(G)$ and $\delta(G)$ the maximum degree and minimum degree of $G$, respectively, %Throughout the paper, 
   and we  simply set $\Delta:=\Delta(G)$ and $\delta:=\delta(G)$ in this paper. 
    %Let $[k]:=\{1, \dots, k\}$ for any positive integer $k$.	
    For any positive integer $k$, let $[k]:=\{1, \dots, k\}$, that is, the set of the first $k$ positive integers.
	A (proper) {\bf $k$-edge-coloring} of a graph $G$ is a map $\phiv:E(G)\to [k]$ that assigns to every edge $e$ of $G$ a color $\phiv(e)\in [k]$ such that no two adjacent edges of $G$ receive the same color.
	Denote by $\CC^k(G)$ the set of all $k$-edge-colorings of $G$.
	The {\bf chromatic index} $\chi'(G)$ is the least integer $k$ such that $\CC^k(G)\ne\emptyset$.
	A graph	$G$ is {\bf edge-chromatic critical} if $\chi'(H)<\chi'(G)$ for every proper subgraph $H$ of $G$. Furthermore, we call $G$ a {\bf $\D$-critical} graph if $G$ is edge-chromatic critical and $\chi'(G)=\D(G)+1$, where $\D =\D(G)$.

	In the 1960s, Vizing \cite{Vizing1965-1} and, independently, Gupta \cite{Gupta1967} proved that $\Delta(G)\le \chi'(G)\le\Delta(G)+1$, which leads to a natural classification of graphs.
	Following Fiorini and Wilson \cite{FW1977}, if $\chi'(G)=\Delta(G)$, then $G$ is of {\bf Class 1}; otherwise, it is of {\bf Class 2}. If a graph $G$ satisfies   $\left|E(G)\right|> \Delta(G)\lfloor \frac{|V(G)|}{2}\rfloor$, then we need at least $\D(G)+1$ colors to get a proper edge coloring of $G$.
	Such graphs are called {\bf overfull}.
	Holyer \cite{Hoyler1981} proved that it is NP-complete to determine whether an arbitrary graph is of Class $1$.
	%The converse is not true in general; nonetheless, a central problem in edge-coloring theory is to identify sufficient conditions that ensure a Class 2 graph is overfull.
    However, a conjecture proposed by Chetwynd and Hilton in 1986 would imply that the chromatic index of graphs $G$ with $\D(G) > \frac{|V(G)|}{3}$ can be determined in polynomial time. The conjecture concerns the overfullness of graphs. In fact, a number of longstanding conjectures listed in \textit{Twenty Pretty Edge Coloring Conjectures} in \cite{SSTF2012} concern deciding when a graph is overfull. The conjecture proposed by Chetwynd and Hilton \cite{CHETWYND1986,CHETWYND1989} in 1986 is stated as follows.
	%Seymour~\cite{Seymour} applying Edmonds' matching polytope theorem, proved that whether a graph $G$ is overfull can be determined in polynomial time.
	%In, Stiebitz showed the , a set of outstanding open problems in edge-coloring theory that includes the Overfull Conjecture.
	
	\begin{Conjecture}[Overfull Conjecture]
		Let $G$ be a graph of Class $2$ with $\D(G)>\frac{|V(G)|}{3}$. Then $G$ contains an overfull subgraph $H$ with $\D(H)=\D(G)$.
	\end{Conjecture}
	
	Note that for the graph $P^*$ obtained from the Petersen graph by deleting one vertex, $\chi'(P^*)=4$ and $\D(P^*)=\frac{|V(P^*)|}{3}=3$, but $P^*$ has no overfull subgraphs of maximum degree $3$ ($P^*$ contains odd cycles, which are overfull).
	Thus, the degree condition $\D(G)>\frac{|V(G)|}{3}$ in the conjecture above is best possible.
	Seymour \cite{Seymour1979} showed that whether a graph $G$ contains an overfull subgraph with maximum degree $\D(G)$ can be determined in polynomial time.
	Niessen \cite{Niessen2001} in 2001 showed that for graphs $G$ with $\D(G)>\frac{|V(G)|}{3}$, there are at most three induced overfull subgraphs with maximum degree $\D(G)$, and it is possible to find one in time polynomial in $|V(G)|$.
	Moreover, when $\D(G)\ge\frac{|V(G)|}{2}$, there is at most one induced overfull subgraph with maximum degree $\D(G)$, and it is possible to find one in time linear in $(|V(G)|+|E(G)|)$.
	Hence, if the  Overfull Conjecture is true, then the NP-complete problem of determining the chromatic index becomes polynomial-time solvable for graphs $G$ with $\D(G)>\frac{|V(G)|}{3}$.
	Furthermore, the Overfull Conjecture would also imply several other longstanding conjectures such as the Just Overfull Conjecture \cite{SSTF2012}, the Vertex-splitting Conjecture \cite{HZ1997}, and Vizing's 2-factor Conjecture, Independence Conjecture and Average Degree Conjecture \cite{Vizing1965-2,Vizing1968} when restricted to graphs $G$ with $\D(G)>\frac{|V(G)|}{3}$.
	
	Although the Overfull Conjecture is very important, there are not many direct research results on it.
	Chetwynd and Hilton \cite{CHETWYND1989} proved that the conjecture holds when $\D(G)\ge |V(G)|-3$.
	Recently, Shan \cite{Shan2024-2} confirmed that for any $0<\varepsilon\le\frac{1}{14}$, there exists a positive integer $n_0$ such that if $G$ is a graph with order $|V(G)|\ge n_0$ and $\D(G)\ge(1-\varepsilon)|V(G)|$, then the conjecture holds.
	There are some other results on the conjecture focusing on minimum degree conditions or both minimum degree and maximum degree conditions.
	%subgraphs induced by maximum degree vertices restricted in some way
	In 2004, Plantholt \cite{Plantholt2004} showed that a graph $G$ with even order and minimum degree $\delta(G) \ge \frac{\sqrt{7}|V(G)|}{3}$ is of Class $2$ if and only if $G$ has an overfull subgraph $H$ with $\Delta(H)=\Delta(G)$. In 2022, Plantholt \cite{Plantholt2022} also obtained the result for large even-order graphs $G$ with  $\delta(G) \ge \frac{2|V(G)|}{3}$. Recently, Plantholt and Shan \cite{Plantholt-Shan2023} showed that for large even-order graphs $G$ whose minimum degree is arbitrarily close to $\frac{1}{2}|V(G)|$ from above,  the Overfull Conjecture holds;  this result has been extended by Shan \cite{Shan2024-1} to the same class of graphs with odd order.
	Since every Class $2$ graph with maximum degree $\Delta$ has a $\D$-critical subgraph, the Overfull Conjecture is equivalent to stating that {\bf every $\D$-critical graph $G$ with $\D(G)>\frac{|V(G)|}{3}$ is overfull}. Note that $\delta(G)\ge2$ for any $\D$-critical graph $G$.
	In 2022, Cao, Chen, Jing, and Shan \cite{CCJS2022} demonstrated that the conjecture holds for $\D$-critical graphs $G$ of order $n$ when $\D(G)-\frac{7\delta(G)}{4}\ge \frac{3n-17}{4}$, which implies the conditions $\Delta(G)\ge\frac{3n-3}{4}$ and $\delta(G)\le\frac{n+1}{7}$ (given $\D(G)\le n-4$ and $\delta(G)\ge 2$). In this paper, we obtain the result below.
   % , which is an improved bound for $n\ge22$.
	
	\begin{Theorem}  \label{Theorem 1}
		Let $G$ be a  $\D$-critical graph of order $n$. If $\D(G)-\frac{5\delta(G)}{3}\ge\frac{2n-7}{3}$, then $G$ is overfull.
	\end{Theorem}
	
	Under the assumptions $\D(G)\le n-4$ and $\delta(G)\ge 2$, it follows readily that the condition $\D(G)-\frac{5\delta(G)}{3}\ge\frac{2n-7}{3}$ implies $\Delta(G)\ge\frac{2n+3}{3}$ and $\delta(G)\le\frac{n-5}{5}$.
    Beyond establishing the  bound above, our proof in this paper introduces new structural and recoloring techniques that may help facilitate further progress on the Overfull Conjecture and related edge-coloring problems.
  %The proof introduces new techniques that may be useful in further investigations of the Overfull.

	\section{Definitions and preliminary results}
	Let $G$ be a graph. For $e \in E(G)$, $G-e$ represents the graph obtained from $G$ by deleting the edge $e$. An edge $e$ is a {\bf critical edge} of $G$ if $\chi'(G-e)<\chi'(G)$. Clearly, if $G$ is of Class $2$ and every edge is critical, then $G$ is edge-chromatic critical. For $u \in V(G)$, let  $d_G(u)$ denote the degree of $u$  in $G$, and let $N_G(u)$ denote the set of neighbors of $u$ in $G$.
	%A $k$-vertex in $G$ is a vertex with degree $k$. For $u \in V(G)$, a $k$-neighbor of  $u$ is a neighbor of $u$ that is a $k$-vertex in $G$.

	Let $G$ be a graph with an edge $e\in E(G)$, and an edge coloring $\phiv\in \CC^{k}(G-e)$ for some integer $k\ge\Delta$. For a vertex $u\in V(G)$, let $\phiv(u)$ denote the set of colors assigned to edges incident with $u$, and $\overline{\varphi}(u) = [k]\backslash\phiv(u)$, that is, the set of colors not assigned to any edge incident with $u$. We call $\phiv(u)$ the set of colors {\bf present} at $u$ and $\overline{\varphi}(u)$ the set of colors {\bf missing} at $u$. Clearly,
	$|\phiv(u)| + |\overline{\varphi}(u)| =k$ for each vertex $u\in V(G)$.
	For a vertex set $X\subseteq V(G)$, we define $\overline{\varphi}(X)=\bigcup_{u\in X} \overline{\varphi}(u)$.
	The set $X$ is called {\bf elementary} with respect to  $\phiv$ or simply {\bf $\phiv$-elementary} if $\overline{\varphi}(u)\cap\overline{\varphi}(v)=\emptyset$ for any two distinct vertices $u,v\in X$.
	For a color $\alpha\in [k]$, let $E_{\varphi,\alpha}(G)$ denote the set of edges colored with $\alpha$.
	Let $\alpha,\beta\in[k]$ be two distinct colors, and let $H$ be the subgraph induced by $E_{\varphi,\alpha}(G)$ and $E_{\varphi,\beta}(G)$.
	Clearly, every component of $H$ is either a path or an even cycle, each of which is referred to as an {\bf $(\alpha,\beta)$-chain} of $G$.
	If we interchange the colors $\alpha$ and $\beta$ on an $(\alpha,\beta)$-chain $C$, then we obtain a new (proper) $k$-edge-coloring of $G-e$, which is also in $\mathcal{C}^k(G-e)$.
	This operation is called a {\bf Kempe change}.
	Furthermore, we say that a chain $C$ has endvertices $u$ and $v$, or that $u$ and $v$ are {\bf $(\alpha,\beta)$-linked}, if $C$ is a path connecting vertices $u$ and $v$.
	For a vertex $u$ of $G$, we denote by $P_u(\alpha,\beta,\varphi)$ or simply $P_u(\alpha,\beta)$ the unique $(\alpha,\beta)$-chain containing the vertex $u$.
	{\bf For two vertices $u, v\in V(G)$, the two chains $P_u(\alpha,\beta,\varphi)$ and $P_v(\alpha,\beta,\varphi)$ are either identical or disjoint}, which is an important fact in our proofs. 
	More generally,  if an $(\alpha,\beta)$-chain is a path containing two vertices
	$a$ and $b$, then we let $P_{[a,b]}(\alpha, \beta,\varphi)$ be its subchain with endvertices $a$ and $b$.
	The operation of interchanging colors $\alpha$ and $\beta$ on the subchain $P_{[a,b]}(\alpha,
	\beta,\varphi)$ is still called a Kempe change, but the resulting coloring may no longer be a proper edge coloring.
	In particular, we may change the color of one edge in our proofs and adopt the notation $uv:\alpha\rightarrow\beta$ to mean that we recolor the edge $uv$ from $\alpha$ to $\beta$. If $\alpha,\beta \in \overline{\varphi}(u)$, then $P_u(\alpha,\beta,\phiv)$ is defined to be the empty chain; by convention, applying a Kempe change on the empty chain does nothing.
	If $x,y \in P_u(\alpha,\beta,\phiv)$ with $x$ positioned between $u$ and $y$ on the path, then we say that $P_u(\alpha,\beta,\phiv)$ meets $x$ before $y$.
	
	%
	%	\begin{figure}[!ht]\label{F1}
		%		\begin{center}
			%			\centering
			%			\scalebox{0.8}{\includegraphics{1.png}}
			%			\caption{(1) The Kempe change on one $(\alpha,\beta)$-chain
				%$P_u(\alpha,\beta,\varphi)$ or $P_v(\alpha,\beta,\varphi)$; (2) The Kempe change on one
				%subchain $P_{[a,b]}(\alpha, \beta,\varphi)$. (Dashed lines represent missing colors at
				%vertices.)}
			%		\end{center}
		%	\end{figure}	

	\begin{Definition}[Multi-fan]
		Let $G$ be a graph with an edge $e=xy_1$, and an edge coloring $\phiv\in \CC^{k}(G-e)$ for some integer $k\ge\Delta$.
		A multi-fan at $x$ with respect to $e$ and $\phiv$ is a sequence $F=(x,e_1,y_1,\dots,e_p,y_p)$ with $p \ge 1$ consisting of edges $e_1, \dots,e_p$ and vertices $x,y_1, \dots ,y_p$ satisfying the following two conditions:
		\begin{itemize}
			\item [{F1.}] The edges $e_1, \dots, e_p$ are distinct, $e_1=e$, and $e_i =xy_i$ for $1\le i\le p$.
			\item [{F2.}] For every edge $e_i$ with $2 \le i \le p$, there exists a vertex $y_j$ with $1 \le j \le i-1$ such that $\phiv(e_i) \in \overline{\varphi}(y_j)$.
		\end{itemize}
	\end{Definition}

	\begin{Definition}[Kierstead path]
		Let $G$ be a graph with an edge $e=y_0y_1$, and an edge coloring $\phiv\in \CC^{k}(G-e)$ for some integer $k\ge\Delta$.	
		A Kierstead path with respect to $e$ and $\phiv$ is a sequence $K=(y_0,e_1,y_1,\dots,e_p,y_p)$ with $p \ge 1$ consisting of edges $e_1, \dots,e_p$ and vertices $y_0,y_1, \dots ,y_p$ satisfying the following two conditions:
		\begin{itemize}
			\item [{K1.}] The vertices $y_0, \dots, y_p$ are distinct, $e_1=e$, and $e_i=y_{i-1}y_i$ for $1 \le i \le p$.
			\item [{K2.}] For every edge $e_i$ with $2 \le i \le p$, there exists a vertex $y_j$  with $0 \le j \le i-2$ such that $\phiv(e_i) \in \overline{\varphi}(y_j)$.
		\end{itemize}
	\end{Definition}	
	
	It is easy to see that a Kierstead path consisting of at most three vertices is also a multi-fan.

	\begin{Lemma}[Vizing's Adjacency Lemma \cite{Vizing1965-1}] \label{L2.1}
		Let $G$ be a graph of Class $2$ with maximum degree $\D$.
		If $e=xy$ is a critical edge of $G$, then $x$ has at least $\D-d_G(y)+1$ neighbors with maximum degree $\D$ in $V(G)\backslash \{y\}$.	
	\end{Lemma}

	\begin{Lemma}[{\cite[Theorem 2.1]{SSTF2012}}]\label{L2.2}
		Let $G$ be a graph of Class $2$ with a critical edge $e_1=xy_1$.
		Let $\phiv\in \CC^{\D}(G-e_1)$ and $F=(x,e_1,y_1, \dots, e_p,y_p)$ be a multi-fan at $x$ with respect to $e_1$ and $\phiv$.
		Then the following statements hold.
		\begin{itemize}
			\item [\textup{(1)}] $V(F)$ is $\phiv$-elementary;
			\item [\textup{(2)}] If $\alpha \in \overline{\varphi}(x)$ and $\beta \in \overline{\varphi}(y_i)$ for $ 1 \le i \le p$, then there exists an $(\alpha,\beta)$-chain with respect to $\phiv$ having endvertices $x$ and $y_i$.
		\end{itemize}
	\end{Lemma}

	\begin{Lemma}[{\cite[Theorem 3.3]{SSTF2012}}]\label{L2.3}
		Let $G$ be a graph of Class $2$ with a critical edge $e_1=y_0y_1$.
		Let $\phiv\in \CC^\Delta(G-e_1)$ and $K=(y_0,e_1,y_1,e_2,y_2,e_3,y_3)$ be a Kierstead path with respect to $e_1$ and $\phiv$. Then the following statements hold.
		\begin{itemize}
			\item [\textup{(1)}] If $\min\{d_G(y_1),d_G(y_2)\} < \Delta$, then $V(K)$ is $\phiv$-elementary;
			%			\item [{\bf (2)}] If $d_G(y_1) \textless \Delta$, then $V(K)$ is $\phiv$-elementary.
			\item [\textup{(2)}] $\left| \overline{\varphi}(y_3) \cap (\overline{\varphi}(y_0) \cup \overline{\varphi}(y_1))\right| \le 1$.	\end{itemize}
	\end{Lemma}

   % The following result is one property of the Kierstead path with $5$ vertices \cite{CCJS2022}, which was first obtained by its first author in his manuscript.
        
    \begin{Lemma} [\cite{CCJS2022}]\label{L2.6}
		Let $G$ be a graph of Class $2$ with a critical edge $ab\in E(G)$. Let $\phiv\in\CC^\Delta(G-ab)$ and $K=(a,ab,b,bu,u,us,s,st,t)$ be a Kierstead path with respect to $ab$ and $\phiv$. If $\left| \overline{\varphi}(t) \cap (\overline{\varphi}(a) \cup\overline{\varphi}(b))\right| \ge 3$, then $d_G(b)=d_G(u)=\Delta$.
	\end{Lemma}

	\begin{Lemma}[\cite{CCJS2022}]\label{L2.4}
		Let $G$ be a $\D$-critical graph of order $n$, and let $a\in V(G)$.
		If $d_G(a)\le\frac{2\D-n+2}{3}$, then for each $v\in V(G)\backslash\{a\}$, either $d_G(v)\le n-\D+2d_G(a)-6$ or $d_G(v)\ge \D-d_G(a)+1$.  Furthermore, if $d_G(v)\ge \D-d_G(a)+1$, then for any $b\in N_G(a)$ with $d_G(b)=\D$ and $\phiv\in \CC^\Delta(G-ab)$, $\left| \overline{\varphi}(v) \cap (\overline{\varphi}(a) \cup\overline{\varphi}(b))\right| \le 1$.
	\end{Lemma}
	
	\begin{Lemma}[\cite{GS2004}]\label{L2.5}
		Let $G$ be a multigraph of order $n$ and $\phiv\in \CC^k(G)$ for some integer $k \ge \Delta$. Then for any color $\alpha \in [k]$, $|\{v \in V(G): \alpha \in \overline{\varphi}(v)\}|\equiv n\pmod{2}$.
	\end{Lemma}
	
%	The  Lemma \ref{L2.5} follows the argument for Case 2 of Theorem 1.2 in \cite{CCJS2022}. Although the assumptions there differ from ours, the reasoning carries over directly; we present the proof in Section 4 for completeness.

	\begin{Lemma} \label{L2.7}
		Let $G$ be a $\D$-critical graph of order $n$ that is not overfull, and let  $ab \in E(G)$. Assume that $d_G(a)=\delta$ and $d_G(b)=\Delta$.
		If $\frac{2n-7+5\delta}{3}\le \Delta \le n-4$ and there exists at most one vertex $t\in V(G)\backslash\{a\}$ such that $d_G(t)\le n-\D+2\delta-6$, then there exist an edge coloring $\phiv\in \CC^\Delta(G-ab)$ and two distinct vertices $x,y\in V(G)\backslash\{a\}$ such that the following three properties hold$:$
        \[
            \begin{gathered}
	               \D-\delta+1\le d_G(x),d_G(y)<\D;\,\,\,\overline{\phiv}(x)\cap(\overline{\phiv}(a)\cup\overline{\phiv}(b))\ne\emptyset;\,\,\,\overline{\phiv}(y)\cap(\overline{\phiv}(a)\cup\overline{\phiv}(b))\ne\emptyset.
            \end{gathered}
        \]
	\end{Lemma}
	
%Note that although Lemma \ref{L2.7} has different hypotheses from those considered in Case 2 of the proof of Theorem 1.2 in \cite{CCJS2022}, it can be proved Lemma \ref{L2.7}. For completeness, we provide the details in Section 4.
Lemma \ref{L2.7} can be obtained using  arguments similar to those in the proof of Case 2 in Theorem 1.2 of \cite{CCJS2022}; its proof will be provided in detail
%for completeness 
in Section 4.

    A {\it short-branch} $H$ is a graph
	with $V(H)=\{a,b,u,x,y\}$ and $E(H)=\{ab,bu,ux,uy\}$.
	
	%(by changing the uncolor edge).

	\begin{Lemma}[{\cite{Qi-Feng2026}}] \label{L2.8}
		Let $G$ be a graph of Class $2$, $H\subseteq G$ be a short-branch with $V(H)=\{a,b,u,x,y\}$ and $E(H)=\{ab,bu,ux,uy\}$, and let $\phiv\in\CC^\Delta(G-ab)$. If
        \[
            \begin{gathered}
	               K=(a,ab,b,bu,u,ux,x)\;{and}\;K^*=(a,ab,b,bu,u,uy,y)
            \end{gathered}
        \]
		are Kierstead paths with respect to $ab$ and $\phiv$, then  $|\overline{\varphi}(x)\cap(\overline{\varphi}(a)\cup\overline{\varphi}(b))|+|\overline{\varphi}(y)\cap(\overline{\varphi}(a)\cup\overline{\varphi}(b))|\le 1$.
	\end{Lemma}
	
Let $G$ be a $\D$-critical graph with $ab\in E(G)$ and $\phiv\in\CC^\D(G-ab)$. A {\it fork} $H \subseteq G$ with respect to  $\phiv$, as defined in \cite{CC2020}, is a subgraph with $V(H)=\{a,b,u,s_1,s_2,t_1,t_2\}$ and $E(H)=\{ab,bu,us_1,us_2,s_1t_1,s_2t_2\}$ such that $\phiv(bu) \in \overline{\varphi}(a)$, $\phiv(us_1),\phiv(us_2) \in \overline{\varphi}(a) \cup \overline{\varphi}(b)$, $\phiv(s_1t_1) \in (\overline{\varphi}(a) \cup \overline{\varphi}(b)) \cap \overline{\varphi}(t_2)$, and $\phiv(s_2t_2) \in (\overline{\varphi}(a) \cup \overline{\varphi}(b)) \cap \overline{\varphi}(t_1)$.

	\begin{Lemma}[\cite{CC2020}] \label{L2.9}
		Let $G$ be a $\D$-critical graph, $ab\in E(G)$, and $\{u,s_1,s_2,t_1,t_2\}\subseteq V(G)$. If $\D\ge d_G(a)+d_G(t_1)+d_G(t_2)+1$, then for any $\phiv\in\CC^\D(G-ab)$, $G$ does not contain the fork on $\{a,b,u,s_1,s_2,t_1,t_2\}$ with respect to  $\phiv$.	
	\end{Lemma}

\begin{figure}[H]
    \centering
    \begin{minipage}[t]{0.45\textwidth}
        \centering
        \vspace{-10pt}
        \scalebox{0.5}{\includegraphics{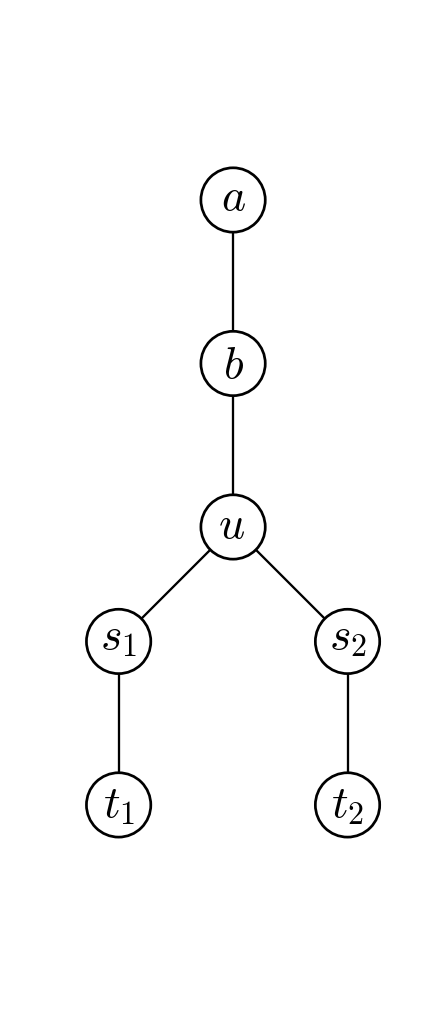}}
        \caption{The branch structure and the two Kierstead paths $K$ and $K^*$ in Lemma \ref{L2.10}.}
        \label{fig:branch}
    \end{minipage}
    \hfill
    \begin{minipage}[t]{0.45\textwidth}
        \centering
        \vspace{-10pt}
        \scalebox{0.5}{\includegraphics{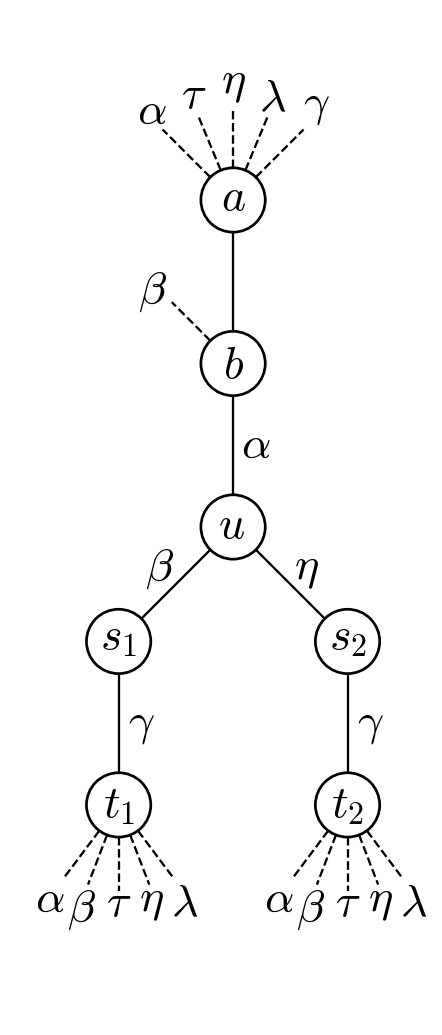}}
        \caption{The edge-coloring configuration on the branch in the proof of Lemma \ref{L2.10}.}
        \label{fig:edge-coloring}
    \end{minipage}
\end{figure}

To improve Lemma~2.8 of \cite{CCJS2022} in the structural sense, we now adjust the structure of the \emph{kite} with  vertex set $\{a,b,c,u,s_1,s_2,t_1,t_2\}$ and  edge set $\{ab,ac,bu,cu,us_1,us_2,s_1t_1,s_2t_2\}$ used there. We refer to the underlying graph of a \emph{fork} 
%$H$ with $V(H)=\{a,b,u,s_1,s_2,t_1,t_2\}$ and $E(H)=\{ab,bu,us_1,us_2,s_1t_1,s_2t_2\}$ 
as a \emph{branch} (see Fig.~\ref{fig:branch}), which omits the auxiliary vertex \(c\) and the edges
\(ac\) and \(cu\). 
%Based on this adjustment, we present the following lemma, and its proof will be given in Section~4.
This structural refinement allows us to
establish the same upper bound on the relevant common missing-color
intersection in a configuration that is more directly applicable to
the proof of Theorem~\ref{Theorem 1}. We give the following lemma and present its proof
in Section~4.

%We improve Lemma~2.8 of \cite{CCJS2022} in the following structural
%sense. The \emph{kite} configuration used there has vertex set
%$\{a,b,c,u,s_1,s_2,t_1,t_2\}$
%and edge set
%$\{ab,ac,bu,cu,us_1,us_2,s_1t_1,s_2t_2\}$.
%We replace this configuration with the following \emph{branch}
%configuration, which omits the auxiliary vertex \(c\) and the edges
%\(ac\) and \(cu\). Formally, a \emph{branch} is the underlying graph
%of a \emph{fork} \(H\) with $V(H)=\{a,b,u,s_1,s_2,t_1,t_2\}$
%and $E(H)=\{ab,bu,us_1,us_2,s_1t_1,s_2t_2\}$
%(Fig.~\ref{fig:branch}). This structural refinement allows us to
%establish the same upper bound on the relevant common missing-color
%intersection in a configuration that is more directly applicable to
%the proof of Theorem~\ref{Theorem 1}. We prove the following lemma
%in Section~4.

	\begin{Lemma}\label{L2.10}
		Let $G$ be a graph of Class $2$, $H\subseteq G$ be a branch with $V(H)=\{a,b,u,s_1,s_2,\\t_1,t_2\}$ and $E(H)=\{ab,bu,us_1,us_2,s_1t_1,s_2t_2\}$, and $\phiv\in\CC^\D(G-ab)$. Assume that
        \[
            \begin{gathered}
	               K=(a,ab,b,bu,u,us_1,s_1,s_1t_1,t_1)\;{and}\;K^*=(a,ab,b,bu,u,us_2,s_2,s_2t_2,t_2)
            \end{gathered}
        \]
		are two Kierstead paths with respect to $ab$ and $\phiv$. If $\phiv(s_1t_1)=\phiv(s_2t_2)$, then $|\overline{\varphi}(t_1)\cap\overline{\varphi}(t_2)\cap(\overline{\varphi}(a)\cup\overline{\varphi}(b))|\le4$.
	\end{Lemma}

	%The proof of Lemma \ref{L2.10} will be given in the Section $4$.
	%It should be noted that the proof of Lemma \ref{L2.8} is also given in our submitted paper \cite{Qi-Yan2025}, and for the sake of completeness, we show it together in this paper.
	
	\section{Proof of Theorem 1}
	
	\begin{proof}
		Since the Overfull Conjecture is true when $\D \ge n-3$, we may assume that $\D \le n-4$.
		Since $\D-\frac{5\delta}{3}\ge\frac{2n-7}{3}$ and $\D \le n-4$, it follows that $\frac{5\delta}{3} \le \Delta-\frac{2n-7}{3} \le (n-4)-\frac{2n-7}{3}=\frac{n-5}{3}$. Hence, $\delta \le \frac{n-5}{5}$.
		Note that $\delta\ge2$ as $G$ is $\D$-critical. Then we have $\D \ge \frac{2n+3}{3}$ and $n\ge15$.
		We select a vertex $a \in V(G)$ with $d_G(a)=\delta$. By Lemma \ref{L2.1}, there exists a vertex $b \in N_G(a)$ with $d_G(b)=\D$. Since $\D-\frac{5\delta}{3}\ge\frac{2n-7}{3}$ and $n \ge 15$,  it can be checked that $\delta< \frac{2\D-n+2}{3}$.
		Thus, by Lemma \ref{L2.4}, for every $v \in V(G) \backslash \{a\}$, we have either $d_G(v) \le n-\D+2\delta-6$ or $d_G(v) \ge \D-\delta+1$.
		%Then we discuss the following two cases to finish the proof.
		
		First, suppose that there exist two distinct vertices $t_1,t_2 \in V(G)\backslash \{a\}$ such that $d_G(t_1),\\d_G(t_2) \le n-\D+2\delta-6$. Notice that $t_i\neq b$ for each $i\in \{1,2\}$. Let $\phiv\in\CC^\Delta(G-ab)$. Since $(a,ab,b)$ is a multi-fan with respect to $ab$ and $\phiv$, $\{a,b\}$ is $\phiv$-elementary by Lemma \ref{L2.2}(1). Thus, $|[\D] \backslash (\overline{\varphi}(a) \cup \overline{\varphi}(b))|=\Delta-|\overline{\varphi}(a) \cup \overline{\varphi}(b)|=\Delta-|\overline{\varphi}(a)|-|\overline{\varphi}(b)|=\Delta-(\Delta-\delta+1)-1=\delta-2$.
		Since $d_G(t_1),d_G(t_2) \ge \delta$ and $|[\D] \backslash(\overline{\varphi}(a) \cup \overline{\varphi}(b))|=\delta-2$, both  $t_1$ and $t_2$ have at least two incident edges colored with colors from $\overline{\varphi}(a) \cup \overline{\varphi}(b)$. Hence, there exist two distinct vertices $s_1 \in N_G(t_1)$ and $s_2 \in N_G(t_2)$ such that $\phiv(s_1t_1),\phiv(s_2t_2) \in \overline{\varphi}(a) \cup \overline{\varphi}(b)$. Since $G$ is $\Delta$-critical, $d_G(s_i) \ge \D-d_G(t_i)+1+1\ge \D-(n-\D+2\delta-6)+1+1 > n-\D+2\delta-6$ for each $i \in \{1,2\}$ by Lemma \ref{L2.1}. Hence, $s_i \neq t_{3-i}$ for each $i\in \{1,2\}$. For each $i\in\{1,2\}$, since $d_G(s_i) > n-\D+2\delta-6\ge \delta=d_G(a)$, we have $s_i \ne a$.
		It follows that $d_G(s_i) \ge \D-\delta+1$ for each $i \in \{1,2\}$ by Lemma \ref{L2.4}.    Note that $|(\overline{\varphi}(a) \cup \overline{\varphi}(b)) \cap \overline{\varphi}(t_i)| \ge 3$ since $|\overline{\varphi}(a) \cup \overline{\varphi}(b)|=\Delta-\delta+2$ and $|\overline{\varphi}(t_i)|=\Delta-d_G(t_i)\ge \Delta-(n-\D+2\delta-6)\ge\Delta-(\D-\delta-1)=\delta+1$ for each $i \in \{1,2\}$.
        We claim that $s_i\neq b$ for each $i \in \{1,2\}$. Otherwise, without loss of generality, suppose that $s_1=b$. Since $\phiv(bt_1) \in \overline{\varphi}(a)$, $(b,ba,a,bt_1,t_1)$ is a multi-fan  with respect to $ab$ and $\phiv$. By Lemma \ref{L2.2}(1), $\{a,b,t_1\}$ is $\phiv$-elementary, contradicting the fact  $|(\overline{\varphi}(a) \cup \overline{\varphi}(b)) \cap \overline{\varphi}(t_1)| \ge 3$.
		
		Between $\{b,s_1,s_2\}$ and $N_G(b)\cap N_G(s_1)\cap N_G(s_2)$, the number of edges colored with  colors from $[\D]\backslash (\overline{\varphi}(a) \cup \overline{\varphi}(b))$ is at most $3(\delta-2)$.
		Note that $|N_G(b)\cap N_G(s_1)\cap N_G(s_2)|=|N_G(b)|+|N_G(s_1)\cap N_G(s_2)|-|N_G(b)\cup (N_G(s_1)\cap N_G(s_2))|\ge \Delta+(|N_G(s_1)|+|N_G(s_2)|-|N_G(s_1)\cup N_G(s_2)|)-n\ge\D+[2(\D-\delta+1)-n]-n=3\D-2n-2\delta+2 \ge 3\delta-5$.
		Therefore, we can find a vertex $u \in N_G(b) \cap N_G(s_1) \cap N_G(s_2)$ such that none of the edges $bu, us_1, us_2$ is colored with a color from $[\Delta] \backslash (\overline{\phiv}(a) \cup \overline{\phiv}(b))$. Note that $u\neq a$; otherwise, $\phiv(as_1),\phiv(as_2) \in \overline{\varphi}(b)$. Since $s_1 \neq s_2$ and $\varphi$ is a proper edge coloring, $\phiv(as_1)\neq \phiv(as_2)$, contradicting the fact  $|\overline{\varphi}(b)| = 1$.  Therefore, we have $u \neq a$ and $\phiv(bu),\phiv(us_1),\phiv(us_2) \in \overline{\varphi}(a) \cup \overline{\varphi}(b)$.
        We claim that $u \notin \{t_1,t_2\}$. Otherwise, without loss of generality,
		suppose that $u=t_1$.
		Then $(a,ab,b,bt_1,t_1,t_1s_1,s_1)$ is a Kierstead path with respect to $ab$ and $\phiv$ and $\{a,b,s_1,t_1\}$ is not $\varphi$-elementary since $|(\overline{\varphi}(a) \cup \overline{\varphi}(b)) \cap \overline{\varphi}(t_1)| \ge 3$, but $d_G(t_1)< \D$, which  contradicts Lemma \ref{L2.3}(1).

		Consequently, the seven vertices $a,b,u,s_1,s_2,t_1,t_2$ are pairwise distinct.  Now we can find a branch $H$ with $V(H)=\{a,b,u,s_1,s_2,t_1,t_2\}$ and $E(H)=\{ab,bu,us_1,us_2,s_1t_1,s_2t_2\}$. Except for the edge $ab$, all other edges in $E(H)$ are colored with colors from $\overline{\varphi}(a) \cup \overline{\varphi}(b)$. Thus, $K=(a,ab,b,bu,u,us_1,s_1,s_1t_1,t_1)$ and $K^*=(a,ab,b,bu,u,us_2,s_2,s_2t_2,t_2)$ are Kierstead paths  with respect to $ab$ and $\phiv$. 
		Note that $|\overline{\varphi}(t_1) \cap \overline{\varphi}(t_2) \cap (\overline{\varphi}(a) \cup \overline{\varphi}(b))| \ge |\overline{\varphi}(t_1) \cap \overline{\varphi}(t_2)|+|\overline{\varphi}(a) \cup \overline{\varphi}(b)|-\D \ge 2[\D-(n-\D+2\delta-6)]-\D+(\D-\delta+2)-\D = 3\D-2n-5\delta+14 \ge 7 >4$.
		If $\phiv(s_1t_1)=\phiv(s_2t_2)$, then this contradicts Lemma \ref{L2.10}. For the case $\phiv(s_1t_1) \ne \phiv(s_2t_2)$,  we make the following claim.
		
		Claim. If $\phiv(s_1t_1) \ne \phiv(s_2t_2)$, then we may assume that $\phiv(s_1t_1) \in \overline{\varphi}(t_2)$ and $\phiv(s_2t_2) \in \overline{\varphi}(t_1)$.
		
		\begin{proof}
			If  $\phiv(s_1t_1) \in \varphi(t_2)$, then we can choose $s_2' \in N_G(t_2)$ with $\phiv(s_1t_1)=\phiv(s_2't_2)$. By the same arguments as for $s_2$, we have $d_G(s_2') \ge \Delta-\delta+1$ and $s_2' \notin \{a,b,t_1\}$. Moreover, $s_2' \ne s_1$, since otherwise the two adjacent edges $s_1t_1$ and $s_1t_2$ would receive the same color. Then we can still find a new vertex $u' \in N_G(b) \cap N_G(s_1) \cap N_G(s_2')$ such that $\phiv(bu'),\phiv(u's_1),\phiv(u's_2') \in \overline{\varphi}(a) \cup \overline{\varphi}(b)$, and $u' \notin \{a,t_1,t_2\}$. Thus, the seven vertices $a,b,u',s_1,s_2',t_1,t_2$ are pairwise distinct. Now we have a branch $H'$ with $V(H')=\{a,b,u',s_1,s_2',t_1,t_2\}$ and $E(H')=\{ab,bu',u's_1,u's_2',s_1t_1,s_2't_2\}$, and both $K=(a,ab,b,bu',u',u's_1,s_1,s_1t_1,t_1)$ and $K^*=(a,ab,b,\\bu',u',u's_2',s_2',s_2't_2,t_2)$ are Kierstead paths  with respect to $ab$ and $\phiv$.
			 Then the proof follows as in the case where $\phiv(s_1t_1)=\phiv(s_2t_2)$. Therefore, we may  assume that $\phiv(s_1t_1) \in \overline{\varphi}(t_2)$, and by similar arguments we may also assume that $\phiv(s_2t_2) \in \overline{\varphi}(t_1)$, as claimed.
		\end{proof}

		Finally, with the assumption in the claim above, the branch $H$ is a fork, but $d_G(a)+d_G(t_1)+d_G(t_2) \le \delta+2(n-\D+2\delta-6)<\D$, which contradicts Lemma \ref{L2.9}.
		
		Therefore, there exists  at most one vertex $t\in V(G)\backslash\{a\}$ such that $d_G(t)\le n-\D+2\delta-6$.
		Suppose to the contrary that $G$ is not overfull.  By Lemma \ref{L2.7}, there exist an edge coloring $\phiv\in\CC^\Delta(G-ab)$ and two distinct vertices $x,y\in V(G)\backslash\{a\}$ satisfying the following properties: $\D-\delta+1\le d_G(x),d_G(y)<\D; \,\,\, \overline{\phiv}(x)\cap(\overline{\phiv}(a)\cup\overline{\phiv}(b))\ne\emptyset;\,\,\, \overline{\phiv}(y)\cap(\overline{\phiv}(a)\cup\overline{\phiv}(b))\ne\emptyset$. 
		We have $|N_G(b)\cap N_G(x)\cap N_G(y)|\ge|N_G(b)|+|N_G(x)\cap N_G(y)|-n\ge|N_G(b)|+(|N_G(x)|+|N_G(y)|-n)-n\ge\D+2(\D-\delta+1)-2n=3\D-2n-2\delta+2\ge3\delta-5$.
		Note that $|[\D]\backslash(\overline{\phiv}(a)\cup\overline{\phiv}(b))|=\delta-2$.
		It follows that no more than $3(\delta-2)$ edges between $\{b,x,y\}$ and $N_G(b)\cap N_G(x)\cap N_G(y)$ are colored with colors from $[\D]\backslash(\overline{\phiv}(a)\cup\overline{\phiv}(b))$. 
        Thus, we can find a vertex $u^*\in N_G(b)\cap N_G(x)\cap N_G(y)$ such that none of the edges $u^*b, u^*x, u^*y$ is colored with a color from $[\Delta] \backslash (\overline{\phiv}(a) \cup \overline{\phiv}(b))$. Note that $u^*\neq a$; otherwise, $\phiv(ax),\phiv(ay) \in \overline{\varphi}(b)$. Since $x \neq y$ and $\varphi$ is a proper edge coloring, $\phiv(ax)\neq \phiv(ay)$, contradicting the fact $|\overline{\varphi}(b)| = 1$.  Therefore, we have $u^*\neq a$ and $\phiv(u^*b),\phiv(u^*x),\phiv(u^*y)\in \overline{\phiv}(a)\cup\overline{\phiv}(b)$. The five vertices $a,b,u^*,x,y$ are pairwise distinct. 
		We have a short-branch $H^*$ with  $V(H^*)=\{a,b,u^*,x,y\}$ and $E(H^*)=\{ab,bu^*,u^*x,u^*y\}$, and both  $K=(a,ab,b,bu^*,u^*,u^*x,x)$ and $K^*=(a,ab,b,bu^*,u^*,u^*y,y)$ are Kierstead paths with respect to $ab$ and $\phiv$.
		Since $|\overline{\varphi}(x)\cap(\overline{\varphi}(a)\cup\overline{\varphi}(b))|+|\overline{\varphi}(y)\cap(\overline{\varphi}(a)\cup\overline{\varphi}(b))|\ge 2$, it gives a contradiction to Lemma \ref{L2.8}. Thus, $G$ is overfull as desired.
	\end{proof}
	
	\section{Proofs of Lemmas 2.7 and 2.10}

%We first prove Lemma 2.7.
We begin with the proof of Lemma \ref{L2.7}.
	\begin{proof}
    Since $\frac{2n-7+5\delta}{3}\le \Delta \le n-4$, we have $\delta\le\frac{n-5}{5}$. Note that $\delta\ge2$ as $G$ is $\Delta$-critical, and thus  $n\ge 15$. Since $\Delta \ge \frac{2n-7+5\delta}{3}$ and $n\ge 15$, it can be checked that $\delta<\frac{2\D-n+2}{3}$. By Lemma \ref{L2.4}, for every  $v\in V(G)\backslash\{a\}$, we have either  $d_G(v)\le n-\D+2\delta-6$ or $d_G(v)\ge \Delta-\delta+1$.  Let $\phiv\in\CC^\Delta(G-ab)$. Since $(a,ab,b)$ is a multi-fan with respect to $ab$ and $\phiv$, $\{a,b\}$ is $\phiv$-elementary by Lemma \ref{L2.2}(1).
		
		First, assume that $n$ is even.  Note that $V(G-ab)=V(G)$. Applying Lemma \ref{L2.5} to $G-ab$ shows that every color in $[\Delta]$ is missing at an even number of vertices of $G$.  Since $\{a,b\}$ is $\phiv$-elementary, every color in $\overline{\phiv}(a)\cup\overline{\phiv}(b)$ is missing  at some vertex in $V(G)\backslash\{a,b\}$. Let $t$ be a vertex with the smallest degree among all vertices in $V(G)\backslash\{a,b\}$. Since $|\overline{\phiv}(a)\cup\overline{\phiv}(b)|=\Delta-\delta+2$ and $|\overline{\phiv}(t)|\le\Delta-\delta$, there exist at least two distinct colors $\alpha,\beta \in (\overline{\phiv}(a)\cup\overline{\phiv}(b))\backslash\overline{\phiv}(t)$. Since $d_G(v)\ge \Delta-\delta+1$ for every $v\in V(G)\backslash\{a,b,t\}$, $\left| \overline{\varphi}(v) \cap (\overline{\varphi}(a) \cup\overline{\varphi}(b))\right| \le 1$ by Lemma \ref{L2.4}. Thus, there exist  two distinct vertices $x,y \in V(G)\backslash\{a,b,t\}$ such that $\alpha \in \overline{\varphi}(x)$ and $\beta \in \overline{\varphi}(y)$. Now $\overline{\phiv}(x)\cap(\overline{\phiv}(a)\cup\overline{\phiv}(b))\ne\emptyset$ and $\overline{\phiv}(y)\cap(\overline{\phiv}(a)\cup\overline{\phiv}(b))\ne\emptyset$, as desired.
		
		Now assume that $n$ is odd. Note that $V(G-ab)=V(G)$. Applying Lemma \ref{L2.5} to $G-ab$ shows that every color in $[\Delta]$ is missing at an odd number of vertices of $G$. Furthermore,
        there exists a color in $[\Delta]$ that is missing at least three vertices of $G$; otherwise, every color in $[\Delta]$ is missing at exactly one vertex of $G$ and thus $|E(G)|=\frac{ \Delta(n-1)}{2}+1$, contradicting the fact that $G$ is not overfull. We consider the following two subcases. If there exists a vertex $t\in V(G)\backslash\{a\}$ such that $d_G(t)\le n-\D+2\delta-6$, then every color in  $\overline{\varphi}(t) \cap (\overline{\varphi}(a) \cup\overline{\varphi}(b))$ is missing  at some vertex in $V(G)\backslash\{a,b,t\}$. Since $\Delta\ge \frac{2n-7+5\delta}{3}$ and $\delta\le\frac{n-5}{5}$, we have $|\overline{\varphi}(t) \cap (\overline{\varphi}(a) \cup\overline{\varphi}(b))| =|\overline{\varphi}(t)|+|\overline{\varphi}(a) \cup\overline{\varphi}(b)|-|\overline{\varphi}(t)\cup(\overline{\varphi}(a) \cup\overline{\varphi}(b))|\ge [\Delta-(n-\Delta+2\delta-6)]+(\Delta-\delta+2)-\Delta=2\Delta-n-3\delta+8 > \delta+1\ge3$.  Since $d_G(v)\ge \Delta-\delta+1$ for every  $v\in V(G)\backslash\{a,b,t\}$, $\left| \overline{\varphi}(v) \cap (\overline{\varphi}(a) \cup\overline{\varphi}(b))\right| \le 1$ by Lemma \ref{L2.4}. Thus, there exist  two distinct vertices $x,y \in V(G)\backslash\{a,b,t\}$ such that $\overline{\phiv}(x)\cap(\overline{\phiv}(a)\cup\overline{\phiv}(b))\ne\emptyset$ and $\overline{\phiv}(y)\cap(\overline{\phiv}(a)\cup\overline{\phiv}(b))\ne\emptyset$, as desired.
		If there is no vertex $t\in V(G)\backslash\{a\}$ such that $d_G(t)\le n-\D+2\delta-6$, then $d_G(v)\ge \Delta-\delta+1$ for every  $v\in V(G)\backslash\{a\}$. Let $\alpha \in [\Delta]$ be a color that is missing at  least three vertices, say $u_1,u_2,u_3$. Note that $\{a,b\}$ is $\phiv$-elementary. If $\alpha \in \overline{\phiv}(a)\cup\overline{\phiv}(b)$,  then we choose two distinct vertices $x,y \in \{u_1,u_2,u_3\}\backslash \{a,b\}$ such that $\overline{\phiv}(x)\cap(\overline{\phiv}(a)\cup\overline{\phiv}(b))\ne\emptyset$ and $\overline{\phiv}(y)\cap(\overline{\phiv}(a)\cup\overline{\phiv}(b))\ne\emptyset$.
		If $\alpha \notin \overline{\phiv}(a)\cup\overline{\phiv}(b)$, then let $\beta \in \overline{\phiv}(a)\cup\overline{\phiv}(b)$. Without loss of generality, $\beta \in \overline{\phiv}(a)$. Since $\beta \in \overline{\phiv}(a)$ and $\alpha \in \phiv(a)$, the chain $P_a(\alpha,\beta,\phiv)$ exists and  one of its endvertices is $a$.  Among $u_1,u_2,u_3$, at most one vertex is $(\alpha,\beta)$-linked with $a$; let $u_1$ be that vertex if it exists.   Both  $P_{u_2}(\alpha,\beta,\phiv)$ and $P_{u_3}(\alpha,\beta,\phiv)$  are disjoint  from $P_a(\alpha,\beta,\phiv)$. We apply a Kempe change on $P_a(\alpha,\beta,\phiv)$.  Denote the resulting (proper) edge coloring by $\phiv$ again. Now $\alpha \in \overline{\phiv}(a)$. There exist two distinct vertices  $x,y \in \{u_2,u_3\}$ such that $\overline{\phiv}(x)\cap(\overline{\phiv}(a)\cup\overline{\phiv}(b))\ne\emptyset$ and $\overline{\phiv}(y)\cap(\overline{\phiv}(a)\cup\overline{\phiv}(b))\ne\emptyset$, as desired. This completes the proof.
	\end{proof}

We next turn to the proof of Lemma \ref{L2.10}.
%We next prove Lemma 2.10.
  %  We now prove Lemma 2.10.
    
	\begin{proof}
		Let $\Gamma=\overline{\varphi}(t_1) \cap \overline{\varphi}(t_2) \cap (\overline{\varphi}(a) \cup \overline{\varphi}(b))$. Suppose that $|\Gamma| \ge 5$. Note that $ab$ is a critical edge.
		  Since $K$ is a Kierstead path and $\left| \overline{\varphi}(t_1) \cap (\overline{\varphi}(a) \cup\overline{\varphi}(b))\right| > 3$, it follows that $d_G(b)=d_G(u)=\Delta$ by Lemma \ref{L2.6}. Let $\alpha,\beta,\tau,\eta,\lambda \in \Gamma$ be distinct colors and $\phiv(s_1t_1)=\phiv(s_2t_2)=\gamma$. Clearly, $\gamma \notin \{\alpha,\beta,\tau,\eta,\lambda\}$. Throughout this proof, whenever we apply Kempe changes simultaneously on the chains at $t_1$ and $t_2$, we apply a Kempe change on each chain if the two chains are disjoint, and apply it only once if the two chains are identical.

          Claim $1$. We may  assume that $\alpha\in\overline{\varphi}(a)$, $\overline{\varphi}(b)=\{\beta\}$, and $\phiv(bu)=\alpha \in \overline{\varphi}(a) \cap \overline{\varphi}(t_1) \cap \overline{\varphi}(t_2)$.
        
         \begin{proof}
        We may first assume that $\alpha\in\overline{\varphi}(a)$ and $\overline{\varphi}(b)=\{\beta\}$. Otherwise, since $|\overline{\varphi}(b)|=1$, we have either 
        $\beta\in\overline{\varphi}(a)$ and $\overline{\varphi}(b)=\{\alpha\}$, or $\alpha,\beta\in\overline{\varphi}(a)$. In the former case, we interchange the names of $\alpha$ and $\beta$, as desired.  In the latter case, let $\overline{\varphi}(b)=\{\alpha'\}$, where $\alpha'\notin \{\alpha,\beta\}$. Since $(a,ab,b)$ is a multi-fan with respect to $ab$ and $\phiv$, there exists a $(\beta,\alpha')$-chain $P_a(\beta,\alpha',\phiv)=P_b(\beta,\alpha',\phiv)$ by Lemma \ref{L2.2}(2). We observe that  $s_1t_1,s_2t_2\notin P_b(\beta,\alpha',\phiv)$; if $\gamma = \alpha'$, then $s_1t_1\in P_{t_1}(\beta,\alpha',\varphi)$ and $s_2t_2\in P_{t_2}(\beta,\alpha',\varphi)$ (both  $P_{t_1}(\beta,\alpha',\varphi)$ and $P_{t_2}(\beta,\alpha',\varphi)$ are disjoint  from $P_b(\beta,\alpha',\varphi)$), and if $\gamma \neq \alpha'$, then this observation is immediate.
        We apply a Kempe change on $P_b(\beta,\alpha',\phiv)$ and denote the resulting (proper) edge coloring by $\phiv$ again. Now $\alpha,\alpha'\in\overline{\varphi}(a)$ and $\overline{\varphi}(b)=\{\beta\}$, as desired. 
        Although this operation may alter the colors of the edges in $K$ and $K^*$, both $K$ and $K^*$ are still Kierstead paths with respect to $ab$ and $\phiv$, because the new color of every changed edge of the branch is still missing at $a$ or $b$. Since $d_G(b)=d_G(u)=\Delta$, $\overline{\varphi}(b)=\{\beta\}$, and $K$ is a Kierstead path, we have $\gamma \in \overline{\varphi}(a)$. Now $\alpha,\tau,\eta,\lambda,\gamma\in\overline{\varphi}(a)$, $\alpha,\beta,\tau,\eta,\lambda \in \overline{\varphi}(t_1)\cap \overline{\varphi}(t_2)$,  $\overline{\varphi}(b)=\{\beta\}$, and $\phiv(s_1t_1)=\phiv(s_2t_2)=\gamma$.

         Now we may assume that $\phiv(bu)=\alpha \in \overline{\varphi}(a) \cap \overline{\varphi}(t_1) \cap \overline{\varphi}(t_2)$. Otherwise, let $\phiv(bu)=\alpha''\neq \alpha$. If $\alpha'' \in \overline{\varphi}(a) \cap \overline{\varphi}(t_1) \cap \overline{\varphi}(t_2)$, then we interchange the names of $\alpha''$ and $\alpha$, as desired.  Now we discuss the case   $\alpha'' \notin \overline{\varphi}(a) \cap \overline{\varphi}(t_1) \cap \overline{\varphi}(t_2)$. Clearly, $\alpha'' \in \overline{\varphi}(a)$ since $K$ is a Kierstead path, and thus $\alpha''\notin \overline{\varphi}(t_1) \cap \overline{\varphi}(t_2)$. Since $(a,ab,b)$ is a multi-fan with respect to $ab$ and $\phiv$, there exists an $(\alpha,\beta)$-chain $P_a(\alpha,\beta,\phiv)=P_b(\alpha,\beta,\phiv)$ by Lemma \ref{L2.2}(2).  We apply a Kempe change on $P_b(\alpha,\beta,\phiv)$. Denote the resulting edge coloring by $\phiv_1$. Now $\beta \in \overline{\varphi_1}(a)$ and $\overline{\varphi_1}(b)=\{\alpha\}$. By Lemma \ref{L2.2}(2), there exists an $(\alpha,\alpha'')$-chain $P_a(\alpha,\alpha'',\phiv_1)=P_b(\alpha,\alpha'',\phiv_1)$. We observe that $s_1t_1,s_2t_2\notin P_b(\alpha,\alpha'',\phiv_1)$; if $\gamma = \alpha''$, then $s_1t_1\in P_{t_1}(\alpha,\alpha'',\phiv_1)$ and $s_2t_2\in P_{t_2}(\alpha,\alpha'',\phiv_1)$  (both  $P_{t_1}(\alpha,\alpha'',\phiv_1)$ and $P_{t_2}(\alpha,\alpha'',\phiv_1)$  are disjoint  from $P_b(\alpha,\alpha'',\phiv_1)$), and if $\gamma \neq \alpha''$, then this observation is immediate.
          We apply a Kempe change on $P_b(\alpha,\alpha'',\phiv_1)$ and obtain an edge coloring $\phiv_2$. Now $\alpha\in\overline{\varphi_2}(a)$, $\overline{\varphi_2}(b)=\{\alpha''\}$, and $\phiv_2(bu)=\alpha$.  By Lemma \ref{L2.2}(2), there exists an $(\alpha'',\beta)$-chain $P_a(\alpha'',\beta,\phiv_2)=P_b(\alpha'',\beta,\phiv_2)$. We observe that  $s_1t_1,s_2t_2 \notin P_b(\alpha'',\beta,\phiv_2)$; if $\gamma = \alpha''$, then $s_1t_1\in P_{t_1}(\alpha'',\beta,\phiv_2)$ and $s_2t_2\in P_{t_2}(\alpha'',\beta,\phiv_2)$, and if $\gamma \neq \alpha''$, then this observation is immediate. We apply a Kempe change on $P_b(\alpha'',\beta,\phiv_2)$. Denote the resulting (proper) edge coloring by $\phiv$ again. Note that $\alpha,\tau,\eta,\lambda,\gamma,\alpha''\in\overline{\varphi}(a)$, $\alpha,\beta,\tau,\eta,\lambda \in \overline{\varphi}(t_1)\cap \overline{\varphi}(t_2)$, $\overline{\varphi}(b)=\{\beta\}$, $\phiv(bu)=\alpha$, and $\phiv(s_1t_1)=\phiv(s_2t_2)=\gamma$. Now $\phiv(bu)=\alpha \in \overline{\varphi}(a) \cap \overline{\varphi}(t_1) \cap \overline{\varphi}(t_2)$, as desired. Moreover, the  operations above preserve the Kierstead path conditions for $K$ and $K^*$.
        \end{proof}

		Claim $2$. We may  assume that $\phiv(us_1)=\beta \in \overline{\varphi}(b) \cap \overline{\varphi}(t_1) \cap \overline{\varphi}(t_2)$.
		
		\begin{proof}
			Assume that $\phiv(us_1)=\beta' \neq \beta$.
			Clearly, $\beta' \in \overline{\varphi}(a)$ since $K$ is a Kierstead path.
			Since  $(a,ab,b)$ is a multi-fan  with respect to $ab$ and $\phiv$, there exist a $(\beta,\gamma)$-chain $P_a(\beta,\gamma,\phiv)=P_b(\beta,\gamma,\phiv)$ and a $(\beta,\beta')$-chain $P_a(\beta,\beta',\phiv)=P_b(\beta,\beta',\phiv)$  by Lemma \ref{L2.2}(2). 
           % Case $1$.  $\beta' \in \overline{\varphi}(t_1) \cap \overline{\varphi}(t_2)$.
            We first consider the case $\beta' \in \overline{\varphi}(t_1) \cap \overline{\varphi}(t_2)$.
             Note that $s_1t_1 \in P_{t_1}(\beta,\gamma,\phiv)$ and $s_2t_2 \in P_{t_2}(\beta,\gamma,\phiv)$   (both  $P_{t_1}(\beta,\gamma,\phiv)$ and $P_{t_2}(\beta,\gamma,\phiv)$  are disjoint  from $P_b(\beta,\gamma,\phiv)$).
            Applying a Kempe change on $P_b(\beta,\gamma,\phiv)$ gives an edge coloring $\phiv_1$, for which $\beta \in \overline{\varphi_1}(a)$ and $\overline{\varphi_1}(b)=\{\gamma\}$. Now there exists a $(\beta',\gamma)$-chain $P_a(\beta',\gamma,\phiv_1)=P_b(\beta',\gamma,\phiv_1)$  by Lemma \ref{L2.2}(2). Note that $us_1,s_1t_1 \in P_{t_1}(\beta',\gamma,\phiv_1)$ and $s_2t_2 \in P_{t_2}(\beta',\gamma,\phiv_1)$. We apply a Kempe change on $P_b(\beta',\gamma,\phiv_1)$ and obtain an  edge coloring $\phiv_2$ with $\gamma \in \overline{\varphi_2}(a)$ and $\overline{\varphi_2}(b)=\{\beta'\}$.
            We interchange the names of $\beta'$ and $\beta$, and denote the resulting (proper) edge coloring  by $\phiv$  again. Note that $\alpha,\tau,\eta,\lambda,\gamma,\beta'\in \overline{\varphi}(a)$,  $\alpha,\beta,\tau,\eta,\lambda,\beta' \in \overline{\varphi}(t_1)\cap \overline{\varphi}(t_2)$, $\overline{\varphi}(b)=\{\beta\}$, $\phiv(bu)=\alpha$, $\phiv(us_1)=\beta$,  and $\phiv(s_1t_1)=\phiv(s_2t_2)=\gamma$.  
            Now $\varphi(us_1) = \beta \in \overline{\varphi}(b) \cap \overline{\varphi}(t_1) \cap \overline{\varphi}(t_2)$, as desired. 
            
    We next discuss the case  $\beta' \notin \overline{\varphi}(t_1) \cap \overline{\varphi}(t_2)$ and consider the following two subcases.  
            %There exists a $(\beta,\beta')$-chain $P_a(\beta,\beta',\phiv)=P_b(\beta,\beta',\phiv)$ with endvertices $a$ and $b$ by Lemma \ref{L2.2}(2). We consider the following two subcases.

             Case $1$. $us_1 \in P_a(\beta,\beta',\phiv)=P_b(\beta,\beta',\phiv)$.

        We apply Kempe changes on $P_{t_1}(\beta,\beta',\phiv)$ and $P_{t_2}(\beta,\beta',\phiv)$. Denote the new  edge coloring by $\phiv_1$  (see Fig.~\ref{fig:claim2}(a)), which satisfies $\beta' \in \overline{\varphi_1}(t_1) \cap \overline{\varphi_1}(t_2)$. By Lemma \ref{L2.2}(2), there exists a $(\beta,\lambda)$-chain $P_a(\beta,\lambda,\phiv_1)=P_b(\beta,\lambda,\phiv_1)$. Applying Kempe changes on $P_{t_1}(\beta,\lambda,\phiv_1)$ and $P_{t_2}(\beta,\lambda,\phiv_1)$ gives an edge coloring  $\phiv_2$  (see Fig.~\ref{fig:claim2}(b)), for which $\beta \in \overline{\varphi_2}(t_1) \cap \overline{\varphi_2}(t_2)$. By Lemma \ref{L2.2}(2), there exists a $(\beta,\gamma)$-chain $P_a(\beta,\gamma,\phiv_2)=P_b(\beta,\gamma,\phiv_2)$. Note that $s_1t_1 \in P_{t_1}(\beta,\gamma,\phiv_2)$ and $s_2t_2 \in P_{t_2}(\beta,\gamma,\phiv_2)$. Applying a Kempe change on $P_b(\beta,\gamma,\phiv_2)$  gives an edge coloring $\phiv_3$  (see Fig.~\ref{fig:claim2}(c)), for which $\beta \in \overline{\varphi_3}(a)$ and $\overline{\varphi_3}(b)=\{\gamma\}$. There exists a $(\gamma,\beta')$-chain $P_a(\gamma,\beta',\phiv_3)=P_b(\gamma,\beta',\phiv_3)$  by Lemma \ref{L2.2}(2). Note that $us_1,s_1t_1 \in P_{t_1}(\gamma,\beta',\phiv_3)$ and $s_2t_2 \in P_{t_2}(\gamma,\beta',\phiv_3)$. Applying a Kempe change on $P_b(\gamma,\beta',\phiv_3)$ gives an edge coloring  $\phiv_4$  (see Fig.~\ref{fig:claim2}(d)). Now $\gamma \in \overline{\varphi_4}(a)$ and $\overline{\varphi_4}(b)=\{\beta'\}$. We interchange the names of $\beta'$ and $\beta$, and
          denote the resulting edge coloring  by $\phiv_5$  (see Fig.~\ref{fig:claim2}(e)). Now $\beta' \in \overline{\varphi_5}(a)$, $\overline{\varphi_5}(b)=\{\beta\}$, and $\phiv_5(us_1)=\beta$.  We apply Kempe changes on $P_{t_1}(\beta',\lambda,\phiv_5)$ and $P_{t_2}(\beta',\lambda,\phiv_5)$. Denote the resulting (proper) edge coloring  by  $\phiv$ again.
        Note that $\alpha,\tau,\eta,\lambda,\gamma,\beta' \in \overline{\varphi}(a)$, $\alpha,\beta,\tau,\eta,\lambda \in \overline{\varphi}(t_1)\cap \overline{\varphi}(t_2)$, $\overline{\varphi}(b)=\{\beta\}$, $\phiv(bu)=\alpha$, $\phiv(us_1)=\beta$, and $\phiv(s_1t_1)=\phiv(s_2t_2)=\gamma$. Now $\varphi(us_1) = \beta \in \overline{\varphi}(b) \cap \overline{\varphi}(t_1) \cap \overline{\varphi}(t_2)$,  as desired.  
 \begin{figure}[H]
	\begin{center}
		\centering
		\scalebox{0.50}{\includegraphics{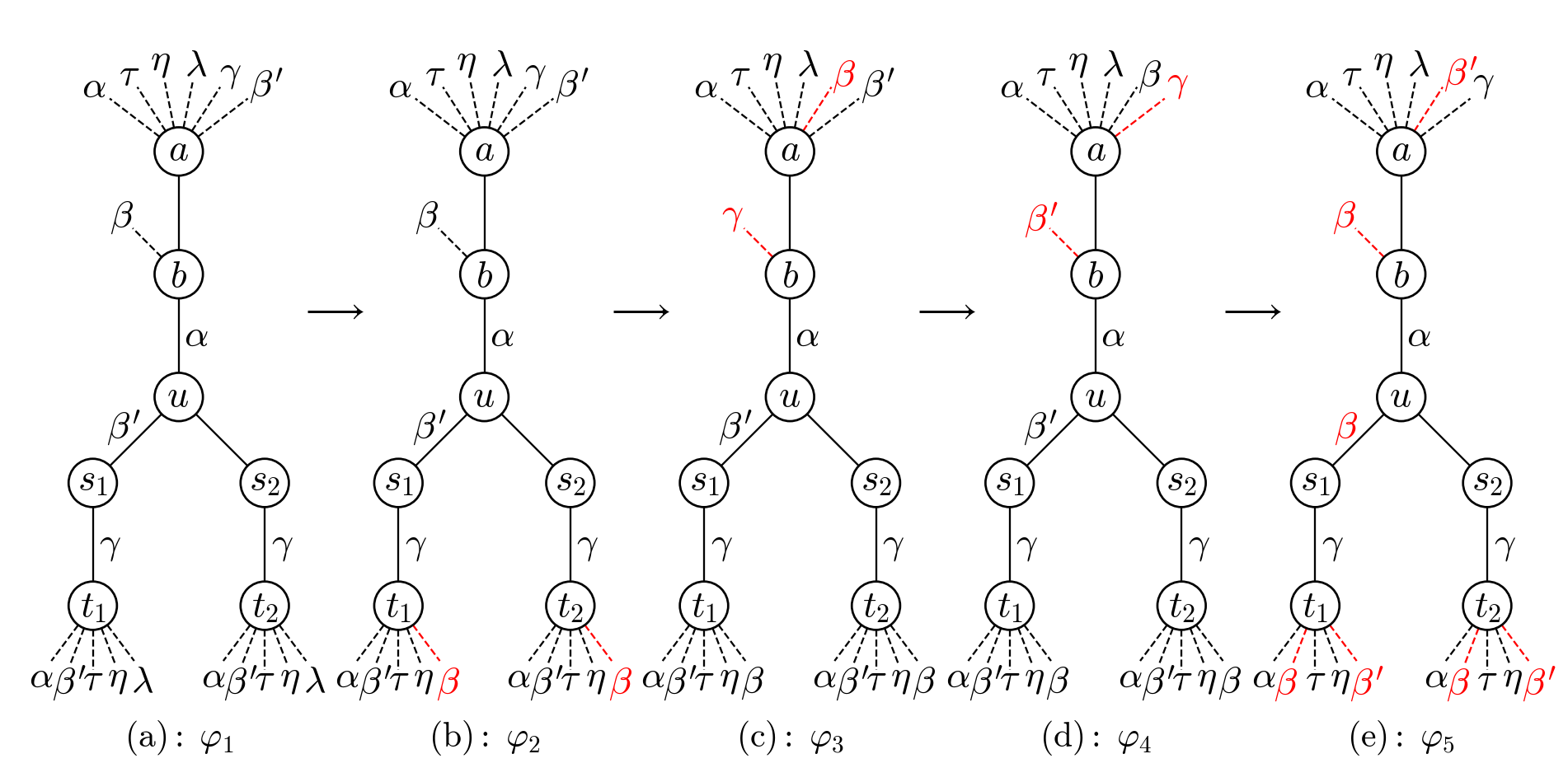}}
		\caption{Selected steps in Case 1 of Claim~2 in the proof of Lemma \ref{L2.10}.}
        \label{fig:claim2}
	\end{center}
\end{figure}

            Case $2$. $us_1 \notin P_a(\beta,\beta',\phiv)=P_b(\beta,\beta',\phiv)$.
            
             We apply a Kempe change on $P_b(\beta,\beta',\phiv)$ and  obtain an edge coloring $\phiv_1$ with $\beta \in \overline{\varphi_1}(a)$ and $\overline{\varphi_1}(b)=\{\beta'\}$. By Lemma \ref{L2.2}(2), there exists an $(\alpha,\beta')$-chain $P_a(\alpha,\beta',\phiv_1)=P_b(\alpha,\beta',\phiv_1)$. Note that $bu,us_1 \in P_b(\alpha,\beta',\phiv_1)$. Apply Kempe changes on $P_{t_1}(\alpha,\beta',\phiv_1)$ and $P_{t_2}(\alpha,\beta',\phiv_1)$, and denote the new edge coloring by $\phiv_2$. Now $\beta' \in \overline{\varphi_2}(t_1) \cap \overline{\varphi_2}(t_2)$ and there exists a $(\gamma,\beta')$-chain $P_a(\gamma,\beta',\phiv_2)=P_b(\gamma,\beta',\phiv_2)$. Note that $us_1,s_1t_1 \in P_{t_1}(\gamma,\beta',\phiv_2)$ and $s_2t_2 \in P_{t_2}(\gamma,\beta',\phiv_2)$. We apply a Kempe change on $P_b(\gamma,\beta',\phiv_2)$ and denote the resulting edge coloring by $\phiv_3$. Now $\beta' \in \overline{\varphi_3}(a)$, $\overline{\varphi_3}(b)=\{\gamma\}$, and there  exists a $(\beta,\gamma)$-chain $P_a(\beta,\gamma,\phiv_3)=P_b(\beta,\gamma,\phiv_3)$.   Applying Kempe changes on $P_{t_1}(\beta,\gamma,\phiv_3)$ and $P_{t_2}(\beta,\gamma,\phiv_3)$ gives an edge coloring $\phiv_4$,  for which  $\phiv_4(s_1t_1)=\phiv_4(s_2t_2)=\beta$ and $\gamma \in \overline{\varphi_4}(t_1) \cap \overline{\varphi_4}(t_2)$. By Lemma \ref{L2.2}(2), there exists an $(\alpha,\gamma)$-chain $P_a(\alpha,\gamma,\phiv_4)=P_b(\alpha,\gamma,\phiv_4)$. Note that $bu \in P_b(\alpha,\gamma,\phiv_4)$. We apply Kempe changes on $P_{t_1}(\alpha,\gamma,\phiv_4)$ and $P_{t_2}(\alpha,\gamma,\phiv_4)$. Denote the new  edge coloring by $\phiv_5$. Now  $\alpha \in \overline{\varphi_5}(t_1) \cap \overline{\varphi_5}(t_2)$ and there exists a $(\lambda,\gamma)$-chain $P_a(\lambda,\gamma,\phiv_5)=P_b(\lambda,\gamma,\phiv_5)$. We apply Kempe changes on $P_{t_1}(\lambda,\gamma,\phiv_5)$ and $P_{t_2}(\lambda,\gamma,\phiv_5)$. Denote the resulting edge coloring  by $\phiv_6$, which satisfies $\gamma \in \overline{\varphi_6}(t_1) \cap \overline{\varphi_6}(t_2)$.  By Lemma \ref{L2.2}(2), there exists a $(\beta,\gamma)$-chain $P_a(\beta,\gamma,\phiv_6)=P_b(\beta,\gamma,\phiv_6)$. Note that $s_1t_1 \in P_{t_1}(\beta,\gamma,\phiv_6)$ and $s_2t_2 \in P_{t_2}(\beta,\gamma,\phiv_6)$. Applying a Kempe change on $P_b(\beta,\gamma,\phiv_6)$  gives an edge coloring $\phiv_7$, for which $\gamma \in \overline{\varphi_7}(a)$ and $\overline{\varphi_7}(b)=\{\beta\}$. There exists a $(\beta,\beta')$-chain $P_a(\beta,\beta',\phiv_7)=P_b(\beta,\beta',\phiv_7)$  by Lemma \ref{L2.2}(2). Note that $us_1,s_1t_1 \in P_{t_1}(\beta,\beta',\phiv_7)$ and $s_2t_2 \in P_{t_2}(\beta,\beta',\phiv_7)$. Applying a Kempe change on $P_b(\beta,\beta',\phiv_7)$ gives an  edge coloring $\phiv_8$. Now $\beta \in \overline{\varphi_8}(a)$ and $\overline{\varphi_8}(b)=\{\beta'\}$. We apply Kempe changes on $P_{t_1}(\beta,\gamma,\phiv_8)$ and $P_{t_2}(\beta,\gamma,\phiv_8)$. Denote the new edge coloring by $\phiv_9$. Now $\phiv_9(s_1t_1)=\phiv_9(s_2t_2)=\gamma$ and $\beta \in \overline{\varphi_9}(t_1)\cap \overline{\varphi_9}(t_2)$.  We interchange the names of $\beta'$ and $\beta$, and denote the resulting edge coloring  by $\phiv_{10}$. Now $\beta' \in \overline{\varphi_{10}}(a)$, $\overline{\varphi_{10}}(b)=\{\beta\}$, and $\phiv_{10}(us_1)=\beta$. We apply Kempe changes on $P_{t_1}(\beta',\lambda,\phiv_{10})$ and $P_{t_2}(\beta',\lambda,\phiv_{10})$.  Denote the resulting (proper) edge coloring by  $\phiv$ again.
            Note that $\alpha,\tau,\eta,\lambda,\gamma,\beta' \in \overline{\varphi}(a)$, $\alpha,\beta,\tau,\eta,\lambda \in \overline{\varphi}(t_1)\cap \overline{\varphi}(t_2)$, $\overline{\varphi}(b)=\{\beta\}$, $\phiv(bu)=\alpha$, $\phiv(us_1)=\beta$, and $\phiv(s_1t_1)=\phiv(s_2t_2)=\gamma$. Now $\varphi(us_1) = \beta \in \overline{\varphi}(b) \cap \overline{\varphi}(t_1) \cap \overline{\varphi}(t_2)$, as desired. 
Moreover, the  operations above preserve the Kierstead path conditions for $K$ and $K^*$.
		\end{proof}

%         \begin{figure}[htbp]
%	\begin{center}
%		\centering
%		\scalebox{0.50}{\includegraphics{claim2.png}}
%		\caption{Partial steps in Case 1 of Claim~2.}
%        \label{fig:claim2}
%	\end{center}
%\end{figure}

		Claim $3$. We may  assume that $\phiv(us_2)=\eta\in\overline{\varphi}(a) \cap \overline{\varphi}(t_1) \cap \overline{\varphi}(t_2)$.
		
		\begin{proof}
			Assume that $\phiv(us_2)=\eta'\neq \eta$. If $\eta'\in\overline{\varphi}(a) \cap \overline{\varphi}(t_1) \cap \overline{\varphi}(t_2)$,  then we interchange the names of $\eta'$ and $\eta$, as desired.  Now we discuss the case  $\eta'\notin\overline{\varphi}(a) \cap \overline{\varphi}(t_1) \cap \overline{\varphi}(t_2)$.
			Clearly, $\eta' \in \overline{\varphi}(a)$ since $K^*$ is a Kierstead path, and thus $\eta' \notin \overline{\varphi}(t_1) \cap \overline{\varphi}(t_2)$. 
			 By Lemma \ref{L2.2}(2), there exists a $(\beta,\gamma)$-chain $P_a(\beta,\gamma,\phiv)=P_b(\beta,\gamma,\phiv)$. Note that $us_1,s_1t_1 \in P_{t_1}(\beta,\gamma,\phiv)$ and $s_2t_2 \in P_{t_2}(\beta,\gamma,\phiv)$.
			Applying a Kempe  change on $P_b(\beta,\gamma,\phiv)$ gives an edge coloring  $\phiv_1$, for which $\beta \in \overline{\varphi_1}(a)$ and $\overline{\varphi_1}(b)=\{\gamma\}$.
			By Lemma \ref{L2.2}(2), there exists a $(\lambda,\gamma)$-chain $P_a(\lambda,\gamma,\phiv_1)=P_b(\lambda,\gamma,\phiv_1)$. Note that $s_1t_1 \in P_{t_1}(\lambda,\gamma,\phiv_1)$ and $s_2t_2 \in P_{t_2}(\lambda,\gamma,\phiv_1)$.
			Applying a Kempe change on $P_b(\lambda,\gamma,\phiv_1)$ gives an edge coloring $\phiv_2$. Now $\gamma \in \overline{\varphi_2}(a)$ and $\overline{\varphi_2}(b)=\{\lambda\}$.
            
			Next we claim that we can obtain a new (proper) edge coloring $\phiv_3$ from $\phiv_2$ such that $\alpha,\beta,\gamma,\tau,\eta,\eta' \in \overline{\varphi_3}(a)$, $\alpha,\beta,\tau,\eta,\eta' \in \overline{\varphi_3}(t_1)\cap \overline{\varphi_3}(t_2)$, $\overline{\varphi_3}(b)=\{\lambda\}$, $\phiv_3(bu)=\alpha$, $\phiv_3(us_1)=\beta$, $\phiv_3(us_2)=\eta'$, and $\phiv_3(s_1t_1)=\phiv_3(s_2t_2)=\gamma$. By Lemma \ref{L2.2}(2), there exists a $(\lambda,\eta')$-chain $P_a(\lambda,\eta',\phiv_2)=P_b(\lambda,\eta',\phiv_2)$.
			We consider the following two subcases to obtain $\phiv_3$.
			If $us_2 \in P_a(\lambda,\eta',\phiv_2)=P_b(\lambda,\eta',\phiv_2)$, then we apply Kempe changes on $P_{t_1}(\lambda,\eta',\phiv_2)$ and $P_{t_2}(\lambda,\eta',\phiv_2)$.
			Denote the resulting (proper)  edge coloring by $\phiv_3$, which satisfies $\eta' \in \overline{\varphi_3}(t_1)\cap\overline{\varphi_3}(t_2)$. Thus, $\varphi_3$ is the desired edge coloring.
			If $us_2\notin P_a(\lambda,\eta',\phiv_2)=P_b(\lambda,\eta',\phiv_2)$, then we apply a Kempe change on $P_b(\lambda,\eta',\phiv_2)$. Denote the resulting  edge coloring by $\phiv_{21}$  (see Fig.~\ref{fig:claim3}(a)), which satisfies $\lambda \in \overline{\varphi_{21}}(a)$ and $\overline{\varphi_{21}}(b)=\{\eta'\}$. By Lemma \ref{L2.2}(2), there exists an $(\alpha,\eta')$-chain $P_a(\alpha,\eta',\phiv_{21})=P_b(\alpha,\eta',\phiv_{21})$. Note that $bu,us_2 \in P_b(\alpha,\eta',\phiv_{21})$.  Applying Kempe changes on $P_{t_1}(\alpha,\eta',\phiv_{21})$ and $P_{t_2}(\alpha,\eta',\phiv_{21})$ gives an edge coloring $\phiv_{22}$  (see Fig.~\ref{fig:claim3}(b)). Now $\eta' \in \overline{\varphi_{22}}(t_1)\cap \overline{\varphi_{22}}(t_2)$ and there exists a $(\gamma,\eta')$-chain $P_a(\gamma,\eta',\phiv_{22})=P_b(\gamma,\eta',\phiv_{22})$. Note that $s_1t_1 \in P_{t_1}(\gamma,\eta',\phiv_{22})$ and $us_2,s_2t_2 \in P_{t_2}(\gamma,\eta',\phiv_{22})$. We apply a Kempe change on $P_b(\gamma,\eta',\phiv_{22})$ and  obtain an  edge coloring $\phiv_{23}$  (see Fig.~\ref{fig:claim3}(c)). Now $\eta' \in \overline{\varphi_{23}}(a)$ and $\overline{\varphi_{23}}(b)=\{\gamma\}$.  By Lemma \ref{L2.2}(2), there exists a $(\lambda,\gamma)$-chain $P_a(\lambda,\gamma,\phiv_{23})=P_b(\lambda,\gamma,\phiv_{23})$. Note that $s_1t_1 \in P_{t_1}(\lambda,\gamma,\phiv_{23})$ and $s_2t_2 \in P_{t_2}(\lambda,\gamma,\phiv_{23})$. We apply a Kempe change on $P_b(\lambda,\gamma,\phiv_{23})$ and  obtain an  edge coloring $\phiv_{24}$  (see Fig.~\ref{fig:claim3}(d)). Now $\gamma \in \overline{\varphi_{24}}(a)$, $\overline{\varphi_{24}}(b)=\{\lambda\}$, and there exists an $(\alpha,\lambda)$-chain $P_a(\alpha,\lambda,\phiv_{24})=P_b(\alpha,\lambda,\phiv_{24})$. Note that $bu \in P_b(\alpha,\lambda,\phiv_{24})$. Applying Kempe changes on $P_{t_1}(\alpha,\lambda,\phiv_{24})$ and $P_{t_2}(\alpha,\lambda,\phiv_{24})$ gives the new (proper) edge coloring $\phiv_3$   (see Fig.~\ref{fig:claim3}(e)), for which $\alpha \in \overline{\varphi_3}(t_1)\cap\overline{\varphi_3}(t_2)$.  Thus, $\varphi_3$ is the desired edge coloring.

			By Lemma \ref{L2.2}(2), there exists a $(\lambda,\tau)$-chain $P_a(\lambda,\tau,\phiv_3)=P_b(\lambda,\tau,\phiv_3)$.
			Applying a Kempe change on $P_b(\lambda,\tau,\phiv_3)$ gives an  edge coloring  $\phiv_4$, for which $\lambda \in \overline{\varphi_4}(a)$ and $\overline{\varphi_4}(b)=\{\tau\}$.  There exists a $(\gamma,\tau)$-chain $P_a(\gamma,\tau,\phiv_4)=P_b(\gamma,\tau,\phiv_4)$ by Lemma \ref{L2.2}(2). Note that $s_1t_1 \in P_{t_1}(\gamma,\tau,\phiv_4)$ and $s_2t_2 \in P_{t_2}(\gamma,\tau,\phiv_4)$. Applying  a Kempe change on $P_b(\gamma,\tau,\phiv_4)$ gives an edge coloring $\phiv_5$. Now $\tau \in \overline{\varphi_5}(a)$, $\overline{\varphi_5}(b)=\{\gamma\}$, and there  exists a $(\beta,\gamma)$-chain $P_a(\beta,\gamma,\phiv_5)=P_b(\beta,\gamma,\phiv_5)$. Note that $us_1,s_1t_1 \in P_{t_1}(\beta,\gamma,\phiv_5)$ and $s_2t_2 \in P_{t_2}(\beta,\gamma,\phiv_5)$. Applying a Kempe change on $P_b(\beta,\gamma,\phiv_5)$  gives an edge coloring $\phiv_6$, for which $\gamma \in \overline{\varphi_6}(a)$ and $\overline{\varphi_6}(b)=\{\beta\}$. We interchange the names of $\eta'$ and $\eta$, and
          denote the new edge coloring  by $\phiv_7$.
           Now $\phiv_7(us_2)=\eta$. We apply Kempe changes on $P_{t_1}(\lambda,\eta',\phiv_7)$ and $P_{t_2}(\lambda,\eta',\phiv_7)$, and denote the new (proper) edge coloring by  $\phiv$ again. Note that $\alpha,\tau,\eta,\lambda,\gamma,\eta'\in \overline{\varphi}(a)$, $\alpha,\beta,\tau,\eta,\lambda \in \overline{\varphi}(t_1)\cap \overline{\varphi}(t_2)$, $\overline{\varphi}(b)=\{\beta\}$,   $\phiv(bu)=\alpha$, $\phiv(us_1)=\beta$, $\phiv(us_2)=\eta$, and $\phiv(s_1t_1)=\phiv(s_2t_2)=\gamma$.
             Now $\phiv(us_2)=\eta\in\overline{\varphi}(a) \cap \overline{\varphi}(t_1) \cap \overline{\varphi}(t_2)$, as desired.
		\end{proof}

         \begin{figure}[htbp]
	\begin{center}
		\centering
		\scalebox{0.50}{\includegraphics{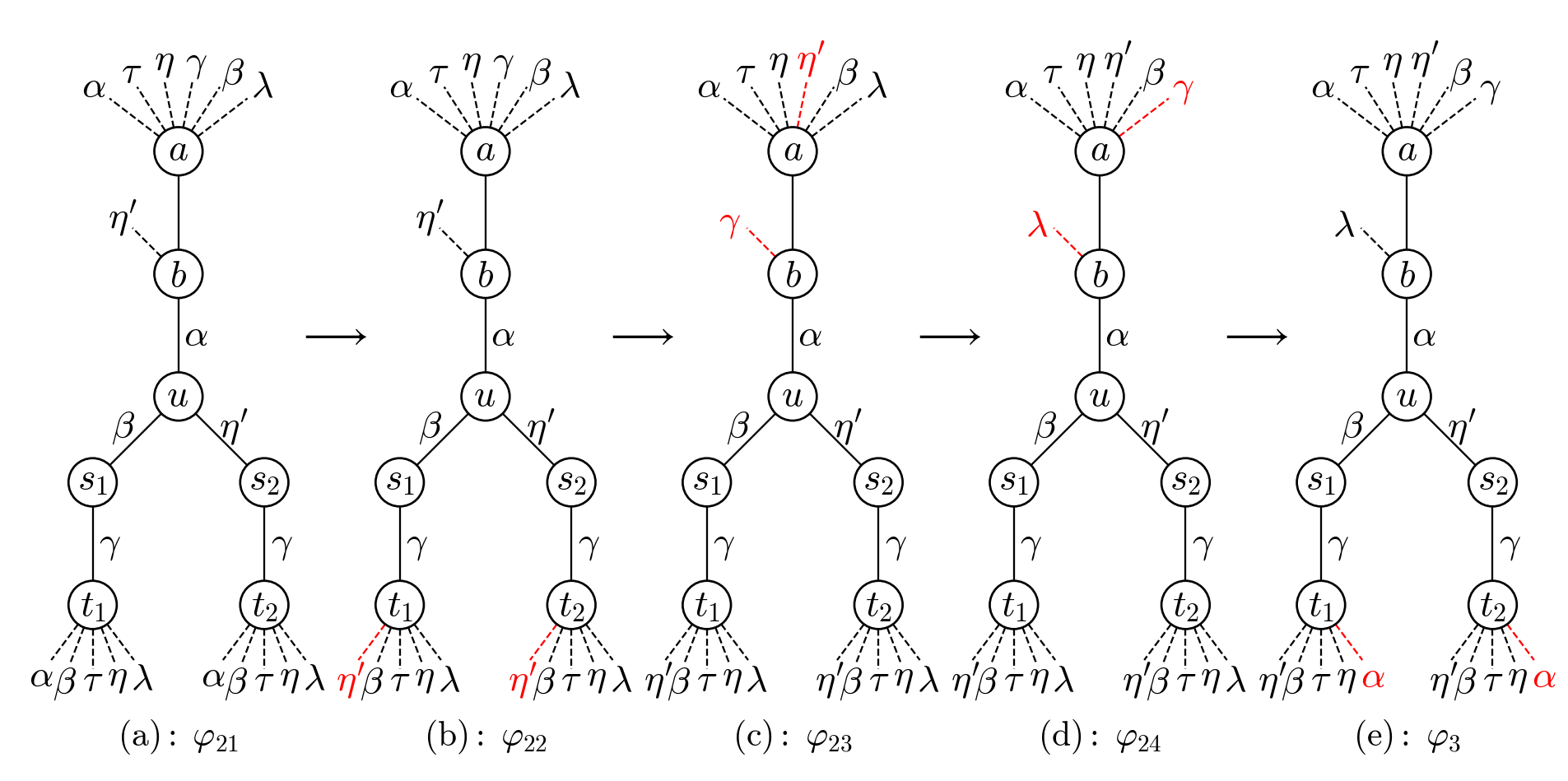}}
		\caption{Selected steps for obtaining $\phiv_3$ in Claim~3 when $us_2\notin P_a(\lambda,\eta',\phiv_2)=P_b(\lambda,\eta',\phiv_2)$ in the proof of Lemma \ref{L2.10}.}
        \label{fig:claim3}
	\end{center}
\end{figure}

        Now $\alpha,\tau,\eta,\lambda,\gamma\in \overline{\varphi}(a)$, $\alpha,\beta,\tau,\eta,\lambda \in \overline{\varphi}(t_1)\cap \overline{\varphi}(t_2)$, $\overline{\varphi}(b)=\{\beta\}$,   $\phiv(bu)=\alpha$, $\phiv(us_1)=\beta$, $\phiv(us_2)=\eta$, and $\phiv(s_1t_1)=\phiv(s_2t_2)=\gamma$  (see Fig.~\ref{fig:edge-coloring}). We present the following two claims using arguments similar to those in the proof of Lemma 2.8  of \cite{CCJS2022}.

		%Notice that $(a,ab,b)$ is a multi-fan with respect to $\varphi$.

        Claim $4$. $bu \in P_{t_1}(\alpha,\gamma,\phiv)$ and $P_{t_1}(\alpha,\gamma,\phiv)$ meets $u$ before $b$.

        \begin{proof}
            We first show  $bu \in P_{t_1}(\alpha,\gamma,\phiv)$. If not,  then  we apply a Kempe change on $P_{t_1}(\alpha,\gamma,\phiv)$  and obtain an edge coloring $\phiv_1$ with $\gamma\in \overline{\varphi_1}(t_1)$ and $\phiv_1(s_1t_1)=\alpha$.  The path $(b,bu,u,us_1,s_1,\\s_1t_1,t_1)$ is an $(\alpha,\beta)$-chain  with endvertices $b$ and $t_1$, contradicting the fact that  $P_b(\alpha,\beta,\phiv_1)=P_a(\alpha,\beta,\phiv_1)$ by Lemma \ref{L2.2}(2).

            Next we claim that $P_{t_1}(\alpha,\gamma,\phiv)$ meets $u$ before $b$. Otherwise, $P_{t_1}(\alpha,\gamma,\phiv)$ meets $b$ before $u$. Let $P_{[s_1,b]}(\alpha,\gamma,\phiv)$ be the subchain of $P_{t_1}(\alpha,\gamma,\phiv)$ with endvertices $s_1$ and $b$. Recoloring  $s_1t_1:\gamma\rightarrow\beta$ gives a (not proper) edge coloring  $\phiv_1$. Now $\gamma \in \overline{\varphi_1}(t_1)$, $\gamma \in \overline{\varphi_1}(s_1)$, and $\phiv_1(s_1t_1)=\phiv_1(us_1)=\beta$. Note that $P_{[s_1,b]}(\alpha,\gamma,\phiv_1)=P_{[s_1,b]}(\alpha,\gamma,\phiv)$. Applying a Kempe change on $P_{[s_1,b]}(\alpha,\gamma,\phiv_1)$ gives a (not proper) edge coloring $\phiv_2$. Now $\alpha \in \overline{\varphi_2}(s_1)$, $\beta, \gamma \in \overline{\varphi_2}(b)$,  and both $bu$ and another edge incident with $b$ are colored with $\alpha$ under $\phiv_2$. Color $ab$ with $\gamma$, recolor $us_1:\beta\rightarrow\alpha$, and recolor $bu:\alpha\rightarrow\beta$.  Denote the resulting  (proper) edge coloring by $\phiv_3$. Under $\phiv_3$, a direct check at all affected vertices $a,b,u,s_1,t_1$ shows that no two incident edges receive the same color. Therefore, $\phiv_3$ is a (proper) $\D$-edge-coloring of $G$,  contradicting  the fact  $\chi'(G)=\D+1$.
        \end{proof}

        Claim $5$. $us_1,us_2 \in P_a(\beta,\eta,\phiv)=P_b(\beta,\eta,\phiv)$.

        \begin{proof}
            Suppose to the contrary that $us_1,us_2 \notin P_a(\beta,\eta,\phiv)=P_b(\beta,\eta,\phiv)$. Let $P$ be the $(\beta,\eta)$-chain containing the edges $us_1$ and $us_2$. We apply a Kempe change on $P$ and obtain an edge coloring $\phiv_1$ with $\phiv_1(us_1)=\eta$ and $\phiv_1(us_2)=\beta$. Since $P_{t_1}(\alpha,\gamma,\phiv_1)=P_{t_1}(\alpha,\gamma,\phiv)$, by Claim 4, we still have $bu \in P_{t_1}(\alpha,\gamma,\phiv_1)$ and $P_{t_1}(\alpha,\gamma,\phiv_1)$  meets $u$ before $b$. Applying an argument similar to that in the proof of Claim 4 to $P_{t_2}(\alpha,\gamma,\phiv_1)$, we obtain that $bu\in P_{t_2}(\alpha,\gamma,\phiv_1)$ and that $P_{t_2}(\alpha,\gamma,\phiv_1)$ meets $u$ before $b$. Thus, the two chains $P_{t_1}(\alpha,\gamma,\phiv_1)$ and $P_{t_2}(\alpha,\gamma,\phiv_1)$ are identical, since both contain $bu$. However, the common chain meets $u$ before $b$ when
            traversed from both $t_1$ and $t_2$, which is impossible since $t_1$ and $t_2$ are its two endvertices.
        \end{proof}

        By Claim 4, let $P_{[t_1,u]}(\alpha,\gamma,\phiv)$ be the subchain of $P_{t_1}(\alpha,\gamma,\phiv)$ with endvertices $t_1$ and $u$. By Claim 5, let $P_{[u,w]}(\beta,\eta,\phiv)$ be the subchain of $P_a(\beta,\eta,\phiv)=P_b(\beta,\eta,\phiv)$ with endvertices $u$ and $w$ such that $us_1\in P_{[u,w]}(\beta,\eta,\phiv)$ but $us_2\notin P_{[u,w]}(\beta,\eta,\phiv)$.  
		%Let $\phiv_1$ be the new edge coloring obtained from $\phiv$ by applying Kempe change on $P_{[u,t]}(\beta,\eta,\phiv)$.
		When $w=a$, we apply a Kempe change on $P_{[u,a]}(\beta,\eta,\phiv)$ and obtain a (not proper) edge coloring $\phiv_{11}$. Now $\beta \in \overline{\varphi_{11}}(u)$, $\beta \in \overline{\varphi_{11}}(a)$, and $\phiv_{11}(us_1)=\phiv_{11}(us_2)=\eta$.  Note that $P_{[t_1,u]}(\alpha,\gamma,\phiv_{11})=P_{[t_1,u]}(\alpha,\gamma,\phiv)$. Applying a Kempe change on $P_{[t_1,u]}(\alpha,\gamma,\phiv_{11})$ gives a (not proper) edge coloring $\phiv_{12}$. Now $\gamma \in \overline{\varphi_{12}}(t_1)$, $\beta,\gamma \in \overline{\varphi_{12}}(u)$, $\phiv_{12}(s_1t_1)=\alpha$, and  both $bu$ and another edge incident with $u$ are colored with $\alpha$ under $\phiv_{12}$. Then color $ab$ with $\alpha$, and recolor $bu:\alpha\rightarrow\beta$, $us_2:\eta\rightarrow\gamma$, and $s_2t_2:\gamma\rightarrow\eta$.   Denote the resulting  (proper) edge coloring by $\phiv_{13}$. Under $\phiv_{13}$, no two incident edges at $a,b,u,s_1,s_2,t_1,t_2$ receive the same color. Therefore, $\phiv_{13}$ is a 
		 (proper) $\D$-edge-coloring of $G$, contradicting the fact $\chi'(G)=\D+1$.
        When $w=b$, we apply a Kempe change on $P_{[u,b]}(\beta,\eta,\phiv)$ and obtain the new (not proper) edge coloring $\phiv_{21}$. Now $\beta \in \overline{\varphi_{21}}(u)$, $\overline{\varphi_{21}}(b)=\{\eta\}$, and $\phiv_{21}(us_1)=\phiv_{21}(us_2)=\eta$. Then color $ab$ with  $\eta$ and uncolor $us_2$. Denote the resulting (proper) edge coloring by $\phiv_{22}$, which satisfies $\overline{\varphi_{22}}(u)=\{\beta\}$ and $\eta \in \overline{\varphi_{22}}(s_2)$. Under $\phiv_{22}$, no two incident edges at $a,b,u,s_1,s_2$ receive the same color.  Since $\phiv_{22}\in\CC^\D(G-us_2)$ and $G$ is of Class 2 with $\chi'(G)=\Delta+1$, $us_2$ is a critical edge.
		Since $(u,us_2,s_2)$ is a multi-fan with respect to $us_2$ and $\phiv_{22}$, there exists a $(\beta,\eta)$-chain $P_u(\beta,\eta,\phiv_{22})=P_{s_2}(\beta,\eta,\phiv_{22})$ by Lemma \ref{L2.2}(2).
		Next we discuss the following two subcases. 
		If $P_u(\beta,\gamma,\phiv_{22})=P_{t_2}(\beta,\gamma,\phiv_{22})$, then we apply a Kempe change on $P_{t_1}(\beta,\gamma,\phiv_{22})$ and obtain a  (proper) edge coloring $\phiv_{23}$ with $\phiv_{23}(s_1t_1)=\beta$ and $\gamma \in \overline{\varphi_{23}}(t_1)$. Hence, the path $(u,us_1,s_1,s_1t_1,t_1)$ is  a $(\beta,\eta)$-chain with endvertices $u$ and $t_1$, contradicting the fact that $P_u(\beta,\eta,\phiv_{23})=P_{s_2}(\beta,\eta,\phiv_{23})$ by Lemma \ref{L2.2}(2).
		If $P_u(\beta,\gamma,\phiv_{22}) \ne P_{t_2}(\beta,\gamma,\phiv_{22})$, then we apply a Kempe change on $P_{t_2}(\beta,\gamma,\phiv_{22})$  and obtain a (proper) edge coloring $\phiv_{24}$ with $\phiv_{24}(s_2t_2)=\beta$ and $\gamma \in \overline{\varphi_{24}}(t_2)$. Color $us_2$ with $\beta$ and recolor $s_2t_2:\beta\rightarrow\eta$. The resulting edge coloring is a (proper) $\D$-edge-coloring of $G$,  contradicting the fact  $\chi'(G)=\D+1$. 
		This completes the proof.
	\end{proof}

	\noindent {\bf Acknowledgements}
	
    The authors are very grateful to the two anonymous reviewers for their valuable comments.
This work was supported by Hebei Natural Science Foundation (A2026208005, A2023208006), Hebei Fund for Introducing Overseas Returnees (C20230357), and Foundation of China Scholarship Council (202508130088).
	
	%	This work was supported by
	%	Hebei Natural Science Foundation A2023208006, Hebei Fund for Introducing Overseas Returnees
	%	C20230357, and NSFC grant 11801135.


\begin{thebibliography}{99}
    
       % \bibitem{C2020}
	%	Y. Cao, {\it Kierstead path of 5 vertices}, manuscript, Georgia State University, 2017.
		


        
		\bibitem{CC2020}
		Y. Cao, G. Chen, {\it On the average degree of edge chromatic critical graphs II}, J. Combin. Theory Ser. B {\bf 145} (2020), 470--486.
		
		\bibitem{CCJS2022}
		Y. Cao, G. Chen, G. Jing, S. Shan, {\it The overfullness of graphs with small minimum degree and large maximum degree}, SIAM J. Discrete Math. {\bf 36} (3) (2022), 2258--2270.
		
		
		
		
		
		\bibitem{CHETWYND1986}
		A.G. Chetwynd, A.J.W. Hilton, {\it Star multigraphs with three vertices of maximum degree}, Math. Proc. Cambridge Philos. Soc. {\bf 100} (2) (1986), 303--317.
		
		\bibitem{CHETWYND1989}
		A.G. Chetwynd, A.J.W. Hilton, {\it The edge-chromatic class of graphs with maximum degree at least $|V|-3$}, Graph theory in memory of G. A. Dirac (Sandbjerg, 1985), Ann. Discrete Math., vol. 41, North-Holland, Amsterdam, 1988, pp. 91--110.
		
		\bibitem{FW1977}
		S. Fiorini, R.J. Wilson, {\it Edge-Colourings of Graphs}, Research Notes in Maths., Pitman, London, 1977.	
		
		\bibitem{GS2004}
		S. Gr\"unewald, E. Steffen, {\it Independent sets and $2$-factors in edge-chromatic-critical graphs}, J. Graph Theory {\bf 45} (2) (2004), 113--118.
		
		\bibitem{Gupta1967}
		R.P. Gupta, {\it Studies in the Theory of Graphs}, PhD thesis, Tata Institute of Fundamental Research, Bombay, 1967.
		
		\bibitem{HZ1997}
		A.J.W. Hilton, C. Zhao, {\it Vertex-splitting and chromatic index critical graphs}, Discrete
		Appl. Math.  {\bf 76} (1997), 205--211.
		
		
		\bibitem{Hoyler1981}
		I. Holyer, {\it The NP-completeness of edge-coloring}, SIAM J. Comput. {\bf 10} (4) (1981), 718--720.
		
		
		\bibitem{Plantholt2004}
		M. Plantholt, {\it Overfull conjecture for graphs with high minimum degree}, J. Graph Theory {\bf 47} (2) (2004),  73--80.
		
		\bibitem{Plantholt2022}
		M. Plantholt, {\it The chromatic index of graphs with large even order n and minimum degree at least $2n/3$}, Discrete Math. {\bf 345} (7) (2022), 112880.
		
		\bibitem{Plantholt-Shan2023}
		M.J. Plantholt, S. Shan,  {\it Edge coloring graphs with large minimum degree}, J. Graph Theory
		{\bf 102} (4) (2023), 611--632.
		
		\bibitem{Qi-Feng2026}
		X. Qi, Y. Feng, {\it On the Hilton-Zhao vertex-splitting conjecture},  Discrete Math. {\bf 349} (2026), 115341.
        %submitted, \url{https://doi.org/10.48550/arXiv.2605.18783}
		
		% \bibitem{Qi-Yan2025}
		%X. Qi, Y. Feng,  {\it On the Hilton-Zhao vertex-splitting conjecture}, submitted, 2025.
		
		\bibitem{Seymour1979}
		P.D. Seymour, {\it On multicolourings of cubic graphs, and conjectures of Fulkerson and Tutte}, Proc.
		London Math. Soc.   {\bf 38} (3) (1979), 423--460.
		
		\bibitem{Shan2024-1}
		S. Shan, {\it The overfull conjecture on graphs of odd order and large minimum degree}, J. Graph Theory
		{\bf106} (2) (2024), 322--351.
		
		\bibitem{Shan2024-2}
		S. Shan, {\it Towards the Overfull Conjecture},  J. Graph Theory {\bf 112} (4) (2026), 370--408.
       
		
		
		
		
		\bibitem{SSTF2012}
		M. Stiebitz, D. Scheide, B. Toft, L.M. Favrholdt, {\it Graph Edge Coloring: Vizing’s Theorem and Goldberg’s Conjecture}, Wiley Series in Discrete Mathematics
and Optimization, John Wiley \& Sons, Inc., Hoboken, NJ, 2012. With a preface by Stiebitz and Toft.

		
		\bibitem{Niessen2001}
		T. Niessen, {\it How to find overfull subgraphs in graphs with large maximum degree. II}, Electron. J.
		Combin. {\bf 8} (2001), \#R7.
		
		
		\bibitem{Vizing1965-1}
		V.G. Vizing, {\it Critical graphs with given chromatic class}, Diskretn. Anal. {\bf 5} (1965), 9--17.
		
		\bibitem{Vizing1965-2}
		V.G. Vizing, {\it The chromatic class of a multigraph}, Kibernetika (Kiev) {\bf 3} (1965), 29--39.
		
		\bibitem{Vizing1968}
		V.G. Vizing, {\it Some unsolved problems in graph theory}, Uspekhi Mat. Nauk. {\bf 23} (6) (1968), 117--134.
		
	\end{thebibliography}
\end{document}